\documentclass[11pt,reqno]{amsart}
\usepackage[a4paper,margin=28mm]{geometry}
\usepackage[T1]{fontenc}
\usepackage[utf8]{inputenc}
\usepackage{lmodern,microtype}
\usepackage{amsmath,amssymb,mathtools}
\usepackage{enumitem,booktabs,array,xcolor,graphicx}
\usepackage[numbers,sort&compress]{natbib}
\usepackage[colorlinks=true,linkcolor=blue!50!black,citecolor=blue!50!black,urlcolor=blue!50!black]{hyperref}
\hypersetup{pdftitle={Optimal adaptive finite element methods in W1p},pdfauthor={Tsogtgerel Gantumur}}
\allowdisplaybreaks
\setlist[enumerate]{label=(\roman*),leftmargin=2.2em}
\newtheorem{theorem}{Theorem}[section]
\newtheorem{proposition}[theorem]{Proposition}
\newtheorem{lemma}[theorem]{Lemma}
\newtheorem{corollary}[theorem]{Corollary}
\theoremstyle{definition}

\newtheorem{assumption}[theorem]{Assumption}
\theoremstyle{remark}
\newtheorem{remark}[theorem]{Remark}

\newcommand{\R}{\mathbb R}
\newcommand{\N}{\mathbb N}
\newcommand{\PP}{\mathbb P}
\newcommand{\adm}{\operatorname{adm}} % configurable admissible-family notation
\newcommand{\conf}{\adm} % compatibility alias for existing source formulas
\newcommand{\jump}[1]{J_{#1}} % face residual, used as \jump{e}(v)
\newcommand{\Id}{\operatorname{Id}}
\newcommand{\ran}{\operatorname{ran}}
\newcommand{\osc}{\operatorname{osc}}

\newcommand{\norm}[1]{\lVert #1\rVert}
\newcommand{\mesh}{T} % triangulations: \mesh, \mesh_0, \mesh'
\newcommand{\faces}[1]{E_{#1}} % interior faces of the mesh in the argument
\newcommand{\errorsymbol}{\mathcal E}
\newcommand{\targeterror}[1]{\errorsymbol_{#1}} % target error (exponent or index list)
\newcommand{\wnorm}[2]{\norm{#1}_{#2}} % W^{1,s}_0 gradient norm: function, exponent

\newcommand{\Sneiberg}{\v{S}ne\u{\i}berg}

\title[Optimal AFEM in $W^{1,p}$]{Optimal adaptive finite element methods in $W^{1,p}$}
\author{Tsogtgerel Gantumur}
\address{Department of Mathematics and Statistics, McGill University, Montr\'eal, Canada}
\email{gantumur.tsogtgerel@mcgill.ca}
\address{Department of Physics, National University of Mongolia, Ulaanbaatar, Mongolia}
\date{October 31, 2026}
\subjclass[2020]{Primary 65N30, 65N50; Secondary 65N15, 46B70}
\keywords{adaptive finite element method, target-norm marking, optimal convergence rates, Ritz stability, residual estimator, general quasi-orthogonality}

\begin{document}

\begin{abstract}
We prove linear convergence and algebraic rate optimality of an adaptive finite element method for the Dirichlet Poisson problem in $W^{1,p}$, $1<p<\infty$. 
We assume well-posedness of the continuous problem and discrete inf-sup stability over the refinement family.
We verify both hypotheses when $p$ is sufficiently close to $2$, without requiring convexity of the domain.
The Galerkin equations remain unchanged, while the residual estimator and D\"orfler marking use the target exponent $p$. 
Linear convergence holds for every fixed marking fraction between zero and one. 
For sufficiently small marking fractions and marked sets of nearly minimal cardinality, the method attains the best algebraic rates for $W^{1,p}$ approximation together with data oscillation. 
The main technical innovation here is a generalization of Feischl's relaxed quasi-orthogonality to a Banach setting, which we prove under the discrete inf-sup stability assumption.
\end{abstract}
\maketitle

\section{Introduction}
\label{sec:intro}

We consider adaptive finite element methods for the Dirichlet problem
\begin{equation}\label{eq:pde}
 -\Delta u=f\quad\hbox{in }\Omega,\qquad u=0\quad\hbox{on }\partial\Omega,
\end{equation}
where $\Omega\subset\R^d$, $d\ge2$, is a bounded polyhedral Lipschitz domain and $f\in L^p(\Omega)$,
for a fixed target exponent $1<p<\infty$.
We assume that the Dirichlet Laplacian is an isomorphism from $W^{1,p}_0(\Omega)$ onto $W^{-1,p}(\Omega)$.
For a conforming triangulation $\mesh$ of $\Omega$, let $u_{\mesh}$ be the standard Galerkin approximation by continuous piecewise polynomials of fixed degree with zero boundary values.
We want to drive mesh refinements by local estimators associated to the target error
\[
 \targeterror{p}(\mesh)
 =\norm{\nabla(u-u_{\mesh})}_{L^p(\Omega)}.
\]
The question is whether this choice leads to a convergent and
rate-optimal method.

In the energy norm, convergence and optimality of adaptive methods are well understood \cite{Doerfler1996,MNS2000,Stevenson2007,CKNS2008}.
The algorithm studied here employs the usual residual based estimators and D\"orfler marking adapted to the target exponent $p$.
We will also follow the proof strategy from the Hilbert setting as far as possible.
This requires upper and lower bounds for the error estimators, together with discrete inf-sup stability in the target Banach norm throughout the admissible mesh family.

Let us mention some earlier precedents on these issues.
Yuan and Li~\cite{YuanLi2011} already establish some of the required
estimator ingredients: reliability and local efficiency up to data
oscillation for residual $W^{1,p}$ estimators.
They consider piecewise affine elements on conforming, shape-regular
triangulations of polygons, with $1<p\le2$ within the continuous
solvability range.
On convex polygons, Badr and Russ~\cite{BadrRuss2014} prove 
discrete $W^{1,p}$ stability for $2<p<2+\varepsilon$, again for
piecewise affine elements on conforming, shape-regular meshes.
For prescribed graded meshes on convex polygons,
Jin, Lee, and Li~\cite{JinLeeLi2017} obtain optimal $W^{1,p}$
approximation rates, while Li~\cite{HLi2017} proves Ritz
stability for $1<p\le\infty$ on sufficiently fine meshes with fixed
corner gradings.
The stability results cited here assume convexity, whereas nonconvex
domains with corner singularities are particularly relevant to
adaptive methods.

In this paper, we establish the additional properties of the
$W^{1,p}$ residual estimators needed for adaptive convergence and
complexity analysis.
Moreover, without assuming convexity of $\Omega$, we prove continuous
and discrete inf-sup stability for all exponents $p$ in an interval
containing $2$, whose size depends on the domain, polynomial degree,
and triangulation class.
Our convergence arguments take continuous and discrete stability
as input, and apply at every exponent for which these properties hold.
The near-two result supplies a nonempty, albeit perturbative, range
of $p$ in which we verify both stability hypotheses.

In the energy norm, Pythagoras' identity bounds the accumulated
Galerkin corrections by the current error.
For $p\ne2$, this identity does not directly provide the required
control in the target norm, but an accumulation estimate allowing
sublinear growth with the number of refinement steps turns out to be sufficient.
Our consideration of such estimates is inspired by the theory of Hilbert space adaptive methods
without energy orthogonality, where they follow from inf-sup stability~\cite{Feischl2022,GantumurGQO}.
Let $(\mesh_j)$ be a nested sequence of uniformly shape-regular
conforming triangulations, and write $u_j=u_{\mesh_j}$ for the
Galerkin solutions.
For all $\ell\ge0$ and $N\ge1$, we prove that
\begin{equation}\label{eq:intro-gqo}
 \sum_{j=\ell}^{\ell+N-1}
 \norm{\nabla(u_{j+1}-u_j)}_{L^p(\Omega)}^p
 \le C N^\alpha \targeterror{p}(\mesh_\ell)^p,
 \qquad 0\le\alpha<1.
\end{equation}
Here we assume a uniform discrete stability at $p$ along the sequence.
The constants $C$ and $\alpha$ may depend on this stability bound,
$p$, the domain, polynomial degree, and mesh regularity, but are
independent of $\ell$ and $N$.
The estimate requires only $u\in W^{1,p}_0(\Omega)$.

We prove \eqref{eq:intro-gqo} by interpolation in Section~\ref{sec:gqo}
and give an independent proof based on uniform convexity in
Appendix~\ref{app:convexity}, that extends the method of \cite{GantumurGQO}.
For the sum in \eqref{eq:intro-gqo}, discrete inf-sup stability alone yields a bound proportional to $N$ (thus $\alpha=1$), whereas energy orthogonality at $p=2$ gives a bound independent of $N$ (thus $\alpha=0$).
Interpolation between these bounds yields $\alpha<1$ at every
exponent strictly between the two endpoints.
To reach the target exponent $p\ne2$, we therefore need stability
at an auxiliary exponent $\hat p$ slightly farther from $2$ than $p$.
Such stability follows from \Sneiberg's theorem, which propagates
discrete stability at $p$ to a neighborhood of $p$,
as we show in Section~\ref{sec:extrapolation}.
Interpolating the bounds at $2$ and $\hat p$ then gives
\eqref{eq:intro-gqo} at the target exponent.

Combining \eqref{eq:intro-gqo} with the residual estimator bounds,
we prove $R$-linear convergence of the estimator and of the target error
combined with data oscillation for every fixed D\"orfler marking
parameter.
Under the mesh assumptions below, we also prove algebraic rate
optimality in the usual sense for sufficiently small marking
parameters and marked sets of nearly minimal cardinality.
Our convergence and complexity arguments follow the standard
approach of~\cite{CFPP2014}.

This paper is organized as follows.
In the next section, we fix the general setting, state the stability
and mesh assumptions, and construct a quasi-interpolation operator.
The stability and accumulation results discussed above are proved
in Sections~\ref{sec:extrapolation}--\ref{sec:gqo}.
Then in Sections~\ref{sec:estimator}--\ref{sec:axioms}, we
establish the required properties of the residual estimators.
The adaptive algorithm and its linear convergence are studied in
Section~\ref{sec:linear}, followed by the proof of rate optimality
in Section~\ref{sec:optimal}.
Finally, Section~\ref{sec:numerics} contains some numerical illustrations.

\section{The continuous problem and its discretization}
\label{sec:setting}

The analysis requires stability of the continuous and discrete
problems in the target norm. We formulate these assumptions,
specify the admissible refinement family, and construct
a quasi-interpolation operator for the stability and residual estimates.

\subsection{Continuous spaces and data}
Let $\Omega\subset\R^d$, $d\ge2$, be a bounded polyhedral Lipschitz domain. For $1<s<\infty$, put
\[
 X_s=W^{1,s}_0(\Omega),\qquad
 \wnorm{v}{s}=\norm{\nabla v}_{L^s(\Omega)},\qquad s'=\frac{s}{s-1}.
\]
Poincar\'e's inequality makes this a Banach norm. Write $X_{s'}^*=W^{-1,s}(\Omega)$, with its operator norm relative to $\wnorm{\cdot}{s'}$. The bilinear form
\begin{equation}\label{eq:bilinear}
 a(v,w)=\int_\Omega\nabla v\cdot\nabla w
\end{equation}
defines a bounded operator $A_s:X_s\to X_{s'}^*$ of norm at most one. Throughout the adaptive analysis $p\in(1,\infty)$ is fixed and $f\in L^p(\Omega)$. 
The functional
\[
 F(w)=\int_\Omega fw
\]
belongs to $X_{p'}^*$ by the H\"older inequality.

\begin{assumption}[Continuous $p$-stability]\label{ass:continuous}
The operator $A_p:X_p\to X_{p'}^*$ is an isomorphism and we have
\begin{equation}\label{eq:continuous-stab}
 \wnorm{v}{p}\le C_{\rm cont}
 \sup_{0\ne w\in X_{p'}}\frac{|a(v,w)|}{\wnorm{w}{p'}}
 \qquad(v\in X_p).
\end{equation}
\end{assumption}

By this assumption, there is a unique $u\in X_p$ with $a(u,w)=F(w)$ for all $w\in X_{p'}$. Surjectivity is required here in addition to the lower inf-sup bound. 
For small $p$ in higher dimensions, $f\in L^p$ need not belong to $H^{-1}$, and the solution is not assumed to lie in $H^1_0$.

\begin{remark}[Continuous stability]\label{rem:continuous-range}
At $p=2$, Assumption~\ref{ass:continuous} follows from the Lax--Milgram theorem. On a bounded Lipschitz domain, Jerison and Kenig~\cite[Theorem~0.5(a)]{JerisonKenig1995} prove that the Dirichlet Laplacian is an isomorphism from $W^{1,p}_0(\Omega)$ onto $W^{-1,p}(\Omega)$ for $p_+'<p<p_+$, where $p_+>4$ if $d=2$ and $p_+>3$ if $d\ge3$. The exponent $p_+$ and the inverse bounds depend on the domain. 
Thus continuous stability is available on a nonempty interval containing $p=2$, without a convexity assumption.
\end{remark}

\subsection{Refinement family}
A conforming simplicial partition $\mesh$ of $\Omega$ is a finite collection of closed $d$-simplices with disjoint interiors whose union is $\bar\Omega$, such that the intersection of two distinct elements is empty or a common boundary subsimplex. Fix a conforming initial partition $\mesh_0$ and a simplicial subdivision rule. Replacing elements by their children, and repeating this operation, defines the family $[\mesh_0]$ of descendant partitions. 
We write $\mesh\preceq\mesh'$ when $\mesh'$ refines $\mesh$ (or equivalently, when $\mesh'\in[\mesh]$). 
Let $\adm(\mesh_0)\subset[\mesh_0]$ be a chosen family of conforming partitions containing $\mesh_0$, called the admissible partitions. Admissibility may impose further restrictions, such as a grading condition. This separation between descendants and admissible partitions follows the setup in \cite{GantumurABEM2013,GantumurMultilevel2017}.

The analysis uses only admissible partitions. For $M\subset\mesh\in\adm(\mesh_0)$, we assume that $\mesh'=\operatorname{refine}(\mesh,M)$ belongs to $\adm(\mesh_0)$, refines $\mesh$, and subdivides every element of $M$ at least once. This operation may consist of subdivision followed by completion to restore admissibility. Thus $M\subset\mesh\setminus\mesh'$. The sets $\mesh\cap\mesh'$, $\mesh\setminus\mesh'$, and $\mesh'\setminus\mesh$ contain the common elements, the refined coarse elements, and the new fine elements, respectively. 
We require the following properties on the set $\adm(\mesh_0)$.
\begin{enumerate}
\item All partitions are uniformly shape regular. The number of elements in each fixed-order vertex patch and the overlap of such patches are uniformly bounded.
\item Every proper descendant $\tau'$ of an element $\tau$ of $\mesh\in\adm(\mesh_0)$ satisfies $|\tau'|\le|\tau|/2$.
\end{enumerate}
For complexity analysis we also assume the following two properties.
\begin{enumerate}[start=3]
\item There is an operation $\oplus:\adm(\mesh_0)\times\adm(\mesh_0)\to\adm(\mesh_0)$. 
Its output $\mesh\oplus\mesh^{\prime}$, called an overlay, refines both arguments ($\mesh\oplus\mesh^{\prime}\in[\mesh]\cap[\mesh^{\prime}]$), 
and satisfies
\begin{equation}\label{eq:overlay}
 \#(\mesh\oplus \mesh')-\#\mesh\le C_{\rm ov}(\#\mesh'-\#\mesh_0).
\end{equation}
\item There is a constant $C_{\rm cls}$ such that, for every
sequence $(\mesh_j)_{j\ge0}$ generated from $\mesh_0$ by
$\mesh_{j+1}=\operatorname{refine}(\mesh_j,M_j)$ with
$M_j\subset\mesh_j$, we have
\begin{equation}\label{eq:closure}
 \#\mesh_\ell-\#\mesh_0
 \le C_{\rm cls}\sum_{j=0}^{\ell-1}\#M_j
 \qquad(\ell\ge1).
\end{equation}
\end{enumerate}

\begin{remark}[Examples of admissible families]\label{rem:admissible-families}
One may take the full conforming newest-vertex bisection family, with compatible initial labels in the relevant dimension \cite{Stevenson2008}. No additional grading restriction is imposed in this case. 
%Corollary~\ref{prop:near-two} supplies a nonempty stability interval about $2$, whose width may depend on the initial mesh, the degree, and the domain. 
Alternatively, admissibility may include a fixed grading condition, with completion preserving that condition. The graded bisection families of \cite[Appendix~A]{DemlowStevenson2011} have the required overlay and closure properties. Their ordinary overlay preserves the grading condition and satisfies \eqref{eq:overlay} with $C_{\rm ov}=1$. 
%For every chosen family, stability at the target exponent remains the separate assumption below. 
\end{remark}

\subsection{Finite element spaces and the Ritz map}
For each integer $m\ge0$, $\PP_m(\sigma)$ denotes the polynomials on $\sigma$ of total degree at most $m$. 
Fix $r\ge1$ and define
\[
 V_{\mesh}=\{v\in C(\bar\Omega):
 v|_\tau\in\PP_r(\tau)\,\,\text{for each}\,\,\tau\in\mesh,\,\,v|_{\partial\Omega}=0\}.
\]
Then let $V_{\mesh,s}$ be $V_{\mesh}$ equipped with the inherited $X_s$ norm. 
We have $V_{\mesh}\subset V_{\mesh'}$ for $\mesh\preceq\mesh'$. 
The Galerkin solution $u_\mesh\in V_\mesh$ and the Ritz map $R_\mesh:X_p\to V_\mesh$ are defined by
\begin{align}
 a(u_{\mesh},w_{\mesh})&=F(w_{\mesh})&& (w_{\mesh}\in V_{\mesh}),\label{eq:galerkin}\\
 a(R_{\mesh}v,w_{\mesh})&=a(v,w_{\mesh})&& (w_{\mesh}\in V_{\mesh},\ v\in X_p).
 \label{eq:ritz}
\end{align}
Note that $u_{\mesh}=R_{\mesh}u$.
The Ritz map is also called the Ritz projector or the
Galerkin projector associated with the bilinear form $a$.
Both definitions are unambiguous by finite-dimensional positive definiteness.

\begin{assumption}[Discrete $p$-stability]\label{ass:discrete}
There is $C_{\rm disc}$, independent of $\mesh\in\conf(\mesh_0)$, such that
\begin{equation}\label{eq:discrete-stab}
 \wnorm{v_{\mesh}}{p}
 \le C_{\rm disc}\sup_{0\ne w_{\mesh}\in V_{\mesh}}
 \frac{|a(v_{\mesh},w_{\mesh})|}{\wnorm{w_{\mesh}}{p'}}
 \qquad(v_{\mesh}\in V_{\mesh}).
\end{equation}
\end{assumption}

This stability assumption gives a uniform bound for the Ritz maps.

\begin{lemma}[Projection stability and compatibility]\label{lem:ritz}
Assumption~\ref{ass:discrete} implies
\begin{equation}\label{eq:ritz-bound}
 K_p:=\sup_{\mesh\in\conf(\mesh_0)}\norm{R_{\mesh}}_{X_p\to X_p}\le C_{\rm disc}.
\end{equation}
Under Assumption~\ref{ass:continuous}, the converse holds with $C_{\rm disc}\le C_{\rm cont}K_p$. 
Moreover, we have
\begin{equation}\label{eq:compatibility}
 R_{\mesh}R_{\mesh'}=R_{\mesh'}R_{\mesh}=R_{\mesh},
\end{equation}
for $\mesh\preceq\mesh'$, and
\begin{equation}\label{eq:cea}
 \targeterror{p}(\mesh)
 \le(1+K_p)\inf_{v_{\mesh}\in V_{\mesh}}\wnorm{u-v_{\mesh}}{p}.
\end{equation}
\end{lemma}

\begin{proof}
For $v\in X_p$, 
using \eqref{eq:discrete-stab},
\eqref{eq:ritz} and H\"older's inequality, we obtain \eqref{eq:ritz-bound} as
\[
 \wnorm{R_{\mesh}v}{p}
 \le C_{\rm disc}\sup_{0\ne w_{\mesh}\in V_{\mesh}}
 \frac{|a(v,w_{\mesh})|}{\wnorm{w_{\mesh}}{p'}}
 \le C_{\rm disc}\wnorm{v}{p}.
\]
For the converse, fix $v_{\mesh}\in V_{\mesh}$ and define $g\in V_{\mesh,p'}^*$ by
$g(w_{\mesh})=a(v_{\mesh},w_{\mesh})$. 
By Hahn--Banach, extend $g$ to a functional in $X_{p'}^*$
without increasing its norm.
Under Assumption~\ref{ass:continuous}, let $z=A_p^{-1}g$.
For every $w_{\mesh}\in V_{\mesh}$, we then have
$a(z,w_{\mesh})=g(w_{\mesh})=a(v_{\mesh},w_{\mesh})$,
so uniqueness in \eqref{eq:ritz} gives $R_{\mesh}z=v_{\mesh}$.
Thus we have
\[
 \wnorm{v_{\mesh}}{p}
 =\wnorm{R_{\mesh}z}{p}
 \le K_p\wnorm{z}{p}
 \le K_p C_{\rm cont}\norm{g}_{V_{\mesh,p'}^*}.
\]
By the definition of the dual norm, this is
\eqref{eq:discrete-stab} with
$C_{\rm disc}\le C_{\rm cont}K_p$.

If $\mesh\preceq\mesh'$, nestedness and reproduction give
$R_{\mesh'}R_{\mesh}=R_{\mesh}$.
For the other order, the Ritz relation on $\mesh'$ gives
$a(R_{\mesh'}v,w_{\mesh})=a(v,w_{\mesh})$ for every
$w_{\mesh}\in V_{\mesh}$. Thus $R_{\mesh}R_{\mesh'}v$
and $R_{\mesh}v$ solve the same discrete problem, proving
\eqref{eq:compatibility}.
Finally, \eqref{eq:cea} is immediate from the fact that
$u-u_{\mesh}=(\Id-R_{\mesh})(u-v_{\mesh})$
for every $v_{\mesh}\in V_{\mesh}$.
\end{proof}

\begin{remark}[Discrete stability]\label{rem:range}
The results recalled in the introduction impose different mesh
restrictions. Badr and Russ~\cite[Theorem~1.1]{BadrRuss2014}
allow conforming mesh families with uniform shape regularity
in their range above two, for piecewise affine elements.
Li~\cite[Theorem~5.5]{HLi2017} treats every $1<p<\infty$
and any fixed polynomial degree on sufficiently fine meshes
with fixed corner gradings.
Both results concern convex polygons.
For Assumption~\ref{ass:discrete}, the bound must hold uniformly
over the entire chosen admissible family, including the meshes
produced by the refinement routine.
Corollary~\ref{prop:near-two} below supplies both stability assumptions
near $p=2$ for the families considered here, without requiring
convexity.
\end{remark}

For $p=2$ and $\mesh\preceq\mesh'$, Galerkin orthogonality
gives Pythagoras' identity
\[
 \wnorm{u-u_{\mesh}}{2}^2
 =\wnorm{u-u_{\mesh'}}{2}^2
  +\wnorm{u_{\mesh'}-u_{\mesh}}{2}^2.
\]
In the standard energy analysis, summing this identity over
successive meshes controls the accumulated Galerkin corrections
by the error at the starting level.
For $p\ne2$, energy orthogonality alone need not even imply
monotonicity in the target norm, as the following algebraic
example shows.

\begin{remark}[Nonmonotonicity of the target error]
\label{ex:nonmonotone}
Let $1<p<\infty$, $p\ne2$, and equip $\R^2$ with the $\ell^p$
norm. Choose $x=(a,b)$ with $ab\ne0$ and $|a|\ne|b|$,
and set $y=(b,-a)$. Let $P_0=0$ and let $P_1$ be the
Euclidean orthogonal projection onto $\operatorname{span}\{y\}$.
Since $x\cdot y=0$, for $u=x+ty$ we have
\[
 u-P_1u=x,\qquad
 \left.\frac{d}{dt}\norm{x+ty}_p^p\right|_{t=0}
 =pab\bigl(|a|^{p-2}-|b|^{p-2}\bigr)\ne0.
\]
For sufficiently small $t$ with sign opposite to this derivative,
we therefore obtain
\[
 \norm{u-P_1u}_p>\norm{u-P_0u}_p.
\]
\end{remark}

Section~\ref{sec:gqo} will establish the sublinear accumulation
estimate needed for the convergence analysis in the target norm.

\subsection{Quasi-interpolation operator}
We will need a local projection onto $V_{\mesh}$ for interpolation
and residual estimation.
The construction uses polynomial moments, as in
\cite{GantumurABEM2013,GantumurMultilevel2017}.
For nested partitions, it leaves functions from the refined space
unchanged outside a fixed neighborhood of the refined elements.

Set $h_\tau=|\tau|^{1/d}$, and let $\faces{\mesh}$ denote
the set of interior faces of $\mesh$ (edges when $d=2$).
For $\tau\in\mesh$ and $k\ge0$, let $\omega_{\mesh}^k(\tau)$
be the interior of the union of all elements reachable from
$\tau$ in at most $k$ steps through elements sharing a vertex.

\begin{lemma}[Local projection and reproduction]\label{lem:local-projector}
There are projectors $I_{\mesh}:X_s\to V_{\mesh,s}$,
defined by the same formula for every $1<s<\infty$,
and a fixed patch order $m$, such that, for every
$v\in X_s$ and $\tau\in\mesh$, we have
\begin{align}
 \norm{\nabla I_{\mesh}v}_{L^s(\tau)}
 &\le C_s\norm{\nabla v}_{L^s(\omega_{\mesh}^{m}(\tau))},
 \label{eq:sz-stab}\\
 \norm{v-I_{\mesh}v}_{L^s(\tau)}
 &\le C_s h_\tau
 \norm{\nabla v}_{L^s(\omega_{\mesh}^{m}(\tau))},
 \label{eq:sz-volume}\\
 \norm{v-I_{\mesh}v}_{L^s(\partial\tau)}
 &\le C_s h_\tau^{1-1/s}
 \norm{\nabla v}_{L^s(\omega_{\mesh}^{m}(\tau))} .
 \label{eq:sz-face}
\end{align}
In particular, we have the global bound
\begin{equation}\label{eq:sz-global}
 \wnorm{I_{\mesh}v}{s}\le C_s\wnorm{v}{s} .
\end{equation}
The constants are uniform in $\mesh$ and locally bounded in $s$.
For every nested pair $\mesh\preceq\mesh'$, there is a set
$Q^+=Q^+(\mesh,\mesh')$ satisfying
\begin{equation}\label{eq:refined-patch}
 Q^+\subset\mesh,\qquad
 \mesh\setminus\mesh'\subset Q^+,\qquad
 \#Q^+\le C_{\rm patch}\#(\mesh\setminus\mesh'),
\end{equation}
such that, for every $\hat v\in V_{\mesh'}$, we have
\begin{equation}\label{eq:local-reproduction}
 (\hat v-I_{\mesh}\hat v)|_\tau=0
 \qquad(\tau\in\mesh\setminus Q^+).
\end{equation}
\end{lemma}

\begin{proof}
Let $\mathcal N_{\mesh}^{\circ}$ be the interior Lagrange
nodes, and let $\phi_z$ be the nodal basis function at $z$.
For each $z\in\mathcal N_{\mesh}^{\circ}$, choose an incident
element $\sigma_z\in\mesh$ and define
$\psi_z\in\PP_r(\sigma_z)$ by
\[
 \int_{\sigma_z}\psi_z q=q(z)
 \qquad\text{for all}\quad q\in\PP_r(\sigma_z).
\]
We then set
\[
 I_{\mesh}v
 =\sum_{z\in\mathcal N_{\mesh}^{\circ}}
   \Big(\int_{\sigma_z}v\psi_z\Big)\phi_z.
\]
These moments define a single linear operator on $L^1(\Omega)$.
For $v\in V_{\mesh}$, each coefficient equals $v(z)$,
and hence $I_{\mesh}v=v$.

Affine scaling gives
$\norm{\psi_z}_{L^{s'}(\sigma_z)}
 \le C_s|\sigma_z|^{-1/s}$.
H\"older's inequality, the bounds for the nodal basis functions,
and a polynomial inverse estimate therefore give
\[
 \norm{I_{\mesh}v}_{L^s(\tau)}
 +h_\tau\norm{\nabla I_{\mesh}v}_{L^s(\tau)}
 \le C_s\norm{v}_{L^s(\omega_{\mesh}^m(\tau))}.
\]
Here $m$ is fixed independently of $s$ and the mesh.
Conformity, shape regularity, and the patch bounds ensure
that the elements involved have diameters comparable to $h_\tau$.

Write $\omega=\omega_{\mesh}^m(\tau)$.
If $\tau$ has no boundary nodes, let $c$ be the average of
$v$ on $\omega$. If $\tau$ has a boundary node, take $c=0$
and choose $m$ large enough that $\omega$ contains an incident
boundary face of diameter comparable to $h_\tau$.
The relevant elements can be joined through a uniformly bounded
number of shared faces. Applying the element Poincar\'e
inequalities along these chains, and using the zero trace of
$v$ on the boundary face in the second case, gives
\[
 \norm{v-c}_{L^s(\omega)}
 \le C_s h_\tau\norm{\nabla v}_{L^s(\omega)}.
\]
In both cases $(I_{\mesh}c)|_\tau=c$.
Applying the preceding local bound to $v-c$ proves
\eqref{eq:sz-stab}, and the triangle inequality gives
\eqref{eq:sz-volume}.
The scaled trace inequality applied to $v-I_{\mesh}v$,
using these two estimates, proves \eqref{eq:sz-face}.
The reference polynomial bounds and the scaled Poincar\'e
and trace constants are uniform when $s$ ranges over a compact
subinterval of $(1,\infty)$.
Summing \eqref{eq:sz-stab} over the mesh and using bounded
patch overlap also gives \eqref{eq:sz-global}.

Finally, define
\(Q^+(\mesh,\mesh')
 =\{\tau\in\mesh:\tau\cap\sigma\ne\varnothing
       \text{ for some }\sigma\in\mesh\setminus\mesh'\}\).
The uniform number of vertex neighbors gives
\eqref{eq:refined-patch}.
If $\tau\notin Q^+$, then $\tau$ and every selected element
$\sigma_z$ associated with its interior nodes are unchanged.
Thus $\hat v|_{\sigma_z}\in\PP_r(\sigma_z)$, and the
corresponding moment equals $\hat v(z)$.
The boundary coefficients agree because $\hat v$ vanishes
on $\partial\Omega$.
Hence $I_{\mesh}\hat v$ and $\hat v$ have identical nodal
values on $\tau$, proving \eqref{eq:local-reproduction}.
\end{proof}

\section{Uniform stability under changes of exponent}
\label{sec:extrapolation}

We prove continuous and discrete stability for exponents near
$p=2$ under the mesh and local projection hypotheses of
Section~\ref{sec:setting}, without requiring convexity of $\Omega$.
More generally, discrete stability at any fixed exponent
persists in a neighborhood of that exponent.
This local openness supplies the neighboring stable exponent
needed to control accumulated Galerkin increments in
Section~\ref{sec:gqo}, without an additional stability assumption.

We use complex-valued Sobolev spaces with the usual Euclidean modulus of the gradient, and complex-linear duals. Extend $a$ complex bilinearly. H\"older's inequality still gives operator norm at most one. A real inverse bound $C_{\rm disc}$ gives a complex inverse bound at most $2C_{\rm disc}$: split a functional on the real subspace into its real and imaginary parts, solve each real problem, and add the two solutions. The same observation applies to the local projectors. All conclusions below restrict to the real spaces.

For a compatible pair of Banach spaces $Y_0,Y_1$, let $[Y_0,Y_1]_\theta$ be their complex interpolation space at $0<\theta<1$. The symbol $\simeq$ between spaces means equality with equivalent norms.

\begin{lemma}[Uniform interpolation scales]\label{lem:interpolation-scales}
Let $1<s_0,s_1<\infty$, 
$1/s=(1-\theta)/s_0+\theta/s_1$ with $0<\theta<1$. Then we have
\begin{align}
 [X_{s_0},X_{s_1}]_\theta&\simeq X_s, 
& [V_{\mesh,s_0},V_{\mesh,s_1}]_\theta&\simeq V_{\mesh,s},
 \label{eq:primal-interpolation}\\
[X_{s_0'}^*,X_{s_1'}^*]_\theta&\simeq X_{s'}^*,
 & [V_{\mesh,s_0'}^*,V_{\mesh,s_1'}^*]_\theta
 &\simeq V_{\mesh,s'}^*.
 \label{eq:dual-interpolation}
\end{align}
For fixed endpoints and $\theta$ in a compact subinterval of
$(0,1)$, the norm equivalence constants are independent of $\mesh\in\adm(\mesh_0)$.
\end{lemma}

\begin{proof}
We first prove the two primal identities.
The continuous identity in \eqref{eq:primal-interpolation}
follows from \cite[Proposition~2.11]{JerisonKenig1995}
at Sobolev order one.
At this order, the supported Bessel-potential spaces in that proposition
are identified with $W^{1,s}_0(\Omega)$ by zero extension.
Poincar\'e's inequality replaces the full Sobolev norms
by our gradient norms.
The proof there gives norm equivalence constants locally
bounded in $\theta$ for fixed endpoints.

For the discrete spaces, let $J_{\mesh}:V_{\mesh,s}\to X_s$
be inclusion and let $I_{\mesh}$ be the common projector from
Lemma~\ref{lem:local-projector}. At the endpoints we have
\[
 I_{\mesh}J_{\mesh}=\Id,\qquad
 \norm{J_{\mesh}}\le1,\qquad
 \norm{I_{\mesh}}_{X_{s_i}\to V_{\mesh,s_i}}\le M_i
 \quad(i=0,1),
\]
where $M_0,M_1$ are independent of $\mesh$.
Interpolation gives bounds of one for $J_{\mesh}$ and
$M_0^{1-\theta}M_1^\theta$ for $I_{\mesh}$ between the
corresponding interpolation spaces.
Since $I_{\mesh}$ reproduces $V_{\mesh}$, these maps identify
$[V_{\mesh,s_0},V_{\mesh,s_1}]_\theta$ with $V_{\mesh}$
equipped with the norm inherited from $[X_{s_0},X_{s_1}]_\theta$.
The continuous identity therefore proves the second identity in
\eqref{eq:primal-interpolation}, uniformly in $\mesh$.

We obtain both dual identities by complex interpolation duality
\cite[Corollary~4.5.2]{BerghLofstrom1976}.
The Sobolev endpoint spaces are reflexive and have the common
dense subspace $C_c^\infty(\Omega)$.
Each discrete pair has the same finite-dimensional underlying
space $V_{\mesh}$, so it also satisfies the density and reflexivity
hypotheses of that result.
The conjugate exponents satisfy
\[
 \frac{1-\theta}{s_0'}+\frac{\theta}{s_1'}
 =1-\frac1s=\frac1{s'}.
\]
Applying duality to the primal identities at $s_0',s_1'$ proves
\eqref{eq:dual-interpolation}.
The constants depend on the ambient Sobolev norm equivalences
and the projector bounds at $s_0,s_1,s_0',s_1'$.
These are independent of $\mesh$ and give the asserted uniformity
in $\theta$.
\end{proof}

We use these uniform interpolation bounds to propagate discrete invertibility from the target exponent to a common neighborhood.
For $1<s<\infty$, let
$A_{\mesh,s}:V_{\mesh,s}\to V_{\mesh,s'}^*$
be defined by $(A_{\mesh,s}v)(w)=a(v,w)$.

\begin{proposition}[Uniform local openness of discrete stability]\label{prop:local-openness}
Suppose Assumption~\ref{ass:discrete} holds at $p\in(1,\infty)$. There are $\varepsilon_p>0$ and $M_p<\infty$ such that
\begin{equation}\label{eq:local-openness}
 \sup_{\mesh\in\conf(\mesh_0)}
 \norm{A_{\mesh,s}^{-1}}_{V_{\mesh,s'}^*\to V_{\mesh,s}}
 \le M_p
 \qquad\text{for}\quad
 |s-p|\le\varepsilon_p .
\end{equation}
The constants depend on $p$, $C_{\rm disc}$, $\Omega$, $r$, and the uniform mesh and local projection bounds. In particular, $K_s:=\sup_{\mesh\in\conf(\mesh_0)}\norm{R_{\mesh}}_{X_s\to X_s}\le M_p$ throughout this neighborhood. No continuous stability at neighboring exponents is required.
\end{proposition}
\begin{proof}
Choose endpoints around the actual target exponent by setting
\[
 \delta_p=\frac12\min\Big\{\frac1p,1-\frac1p\Big\},\qquad
 \frac1{s_0}=\frac1p+\delta_p,\qquad
 \frac1{s_1}=\frac1p-\delta_p.
\]
Then $1<s_0<p<s_1<\infty$, and $s=p$ in
Lemma~\ref{lem:interpolation-scales} when $\theta=\frac12$.
For each $\mesh$, the operators $A_{\mesh,s_0}$ and
$A_{\mesh,s_1}$ are realizations of the same linear map,
with different domain and range norms.
H\"older's inequality gives $\|A_{\mesh,s_i}\|\le1$
for $i=0,1$.
Assumption~\ref{ass:discrete}, after complexification,
gives $\|A_{\mesh,p}^{-1}\|\le2C_{\rm disc}$.

For clarity, write $F_{\mesh,\theta}=[V_{\mesh,s_0},V_{\mesh,s_1}]_\theta$ and $G_{\mesh,\theta}=[V_{\mesh,s_0'}^*,V_{\mesh,s_1'}^*]_\theta$. Lemma~\ref{lem:interpolation-scales} gives positive mesh-independent constants $a_V,b_V,a_G,b_G$ such that
\[
 a_V\wnorm{v}{s}\le\norm v_{F_{\mesh,\theta}}
 \le b_V\wnorm{v}{s},\qquad
 a_G\norm g_{V_{\mesh,s'}^*}\le\norm g_{G_{\mesh,\theta}}
 \le b_G\norm g_{V_{\mesh,s'}^*} ,
\]
for $\frac14\le\theta\le\frac34$.
Thus the inverse norm between interpolation spaces at $\theta=\frac12$ is at most
\[
 B_p:=\max\{1,2b_V C_{\rm disc}/a_G\},
\]
uniformly over all meshes. We use the quantitative form of
\Sneiberg's stability theorem~\cite{Sneiberg1974,CaoSagher1990}
given in \cite[Theorem~A.1]{AuscherBortzEgertSaari2019}.
There are absolute constants $0<c_*\le\frac14$ and $C_*\ge1$
such that the inverse norm is at most $C_*B_p$ whenever
\[
 \big|\theta-\tfrac12\big|\le\rho_p,
 \qquad \rho_p:=\frac{c_*}{1+B_p}.
\]
For example, one may take $c_*=\frac1{24}$ and $C_*=8$.
This interval lies in $[\frac14,\frac34]$. 
The norm equivalences above give $M_p=C_*B_pb_G/a_V$, and one may take $\varepsilon_p=\delta_p\rho_p$ in \eqref{eq:local-openness}. Since the norm equivalence constants are independent of $C_{\rm disc}$,
these choices give $\varepsilon_p\ge c/(1+C_{\rm disc})$
and $M_p\le C(1+C_{\rm disc})$, with the remaining data fixed.
Every quantity used to choose this smaller neighborhood is independent of the mesh. The zero-dimensional spaces cause no difficulty. Restricting to real data proves the assertion, and the Ritz bound follows as in Lemma~\ref{lem:ritz}.
\end{proof}

\begin{corollary}[Uniform stability near $p=2$]\label{prop:near-two}
There are $\varepsilon>0$ and $C<\infty$, depending only on
$\Omega$, $r$, and the admissible family $\adm(\mesh_0)$,
such that Assumptions~\ref{ass:continuous} and~\ref{ass:discrete}
hold for every $|p-2|<\varepsilon$, with
$C_{\rm cont},C_{\rm disc}\le C$.
\end{corollary}

\begin{proof}
At $p=2$, Lax--Milgram gives both stability assumptions with
constant one, on $X_2$ and on every $V_{\mesh,2}$.
Proposition~\ref{prop:local-openness} extends the discrete bound
to nearby exponents.
The same quantitative \Sneiberg\ theorem, applied to the continuous
interpolation scales of Lemma~\ref{lem:interpolation-scales},
extends the continuous isomorphism and its inverse bound.
Taking a smaller common neighborhood proves the assertion.
\end{proof}

\section{Sublinear accumulation of Galerkin increments}
\label{sec:gqo}

Here we prove the accumulation estimate by interpolation.
An independent proof based on uniform convexity is given in
Appendix~\ref{app:convexity}.

\begin{theorem}[Accumulation by interpolation]\label{thm:interpolation-accumulation}
Suppose Assumption~\ref{ass:discrete} holds. There are constants $C_{\rm gqo}>0$ and $0\le\alpha<1$ such that, for every nested sequence $(\mesh_j)_{j\ge0}$ in $\conf(\mesh_0)$, every $N\ge1$, and $v\in X_p$, we have
\begin{equation}\label{eq:interpolated-analysis}
 \sum_{j=1}^N\wnorm{(R_{\mesh_j}-R_{\mesh_{j-1}})v}{p}^p
 \le C_{\rm gqo}N^\alpha\wnorm{v}{p}^p.
\end{equation}
The constants depend on $p$, $C_{\rm disc}$, $\Omega$, $r$, and the uniform mesh and local projection bounds, but not on $N$ or the refinement history.
\end{theorem}
\begin{proof}
For $p=2$, the assertion follows directly from orthogonality.
Suppose henceforth that $p\ne2$, and choose $\hat p$ in the neighborhood of Proposition~\ref{prop:local-openness} so that
\[
 2<p<\hat p\quad\hbox{or}\quad1<\hat p<p<2,
 \qquad
 \frac1p=\frac{1-\theta}{2}+\frac\theta{\hat p},\quad 0<\theta<1.
\]
Let $C_I$ bound the identity map $X_p\to[X_2,X_{\hat p}]_\theta$ on the complex spaces. 
We will prove the estimate \eqref{eq:interpolated-analysis} with
\begin{equation}\label{eq:interpolation-constants}
 \alpha=\frac{\theta p}{\hat p}=\frac{p-2}{\hat p-2}\in(0,1),
 \qquad C_{\rm gqo}=C_I^p(4K_{\hat p})^{\theta p}.
\end{equation}
Here $C_I$ depends on $\Omega,p,\hat p$, whereas $K_{\hat p}$ is controlled by Proposition~\ref{prop:local-openness}.

Write $R_j=R_{\mesh_j}$. Its defining moments against finite element test functions make sense for every $v\in X_s$, $1<s<\infty$. The stiffness matrix is independent of $s$, so these are compatible realizations of one linear map. For fixed $N$, define
\[
 D_Nv
 =\bigl(\nabla(R_j-R_{j-1})v\bigr)_{j=1}^N,
 \qquad Y_{N,s}=L^s(\Omega\times\{1,\ldots,N\};\mathbb C^d),
\]
with counting measure in the index variable. These endpoint operators agree on the intersection and define a single operator on $X_2+X_{\hat p}$.

Since the details are mutually orthogonal in $X_2$, we have
\begin{equation}\label{eq:detail-hilbert-endpoint}
 \norm{D_Nv}_{Y_{N,2}}^2
 =\wnorm{(R_N-R_0)v}{2}^2
 \le\wnorm{v}{2}^2.
\end{equation}
At the other endpoint, $K_{\hat p}$ controls each Ritz map
separately. The triangle inequality bounds each difference by
$2K_{\hat p}\wnorm{v}{\hat p}$, so summing the $\hat p$th
powers gives
$\norm{D_Nv}_{Y_{N,\hat p}}^{\hat p}\le(2K_{\hat p})^{\hat p}N \wnorm{v}{\hat p}^{\hat p}$
for real $v$. Splitting real and imaginary parts bounds the complexified Ritz maps by $2K_{\hat p}$ and hence gives
\begin{equation}\label{eq:detail-stable-endpoint}
 \norm{D_Nv}_{Y_{N,\hat p}}
 \le4K_{\hat p}N^{1/\hat p}\wnorm{v}{\hat p}
 \qquad(v\in X_{\hat p}).
\end{equation}
The extra factor only accounts for complexification. No sharp constant is needed.

Interpolate this same operator using Lemma~\ref{lem:interpolation-scales} and the isometric identity
\[
 [Y_{N,2},Y_{N,\hat p}]_\theta=Y_{N,p}.
\]
The latter is ordinary $L^s$ interpolation with a fixed Euclidean fiber \cite[Theorem~5.1.2]{BerghLofstrom1976}. Its constant is one for every $N$. It follows that
\[
 \norm{D_Nv}_{Y_{N,p}}
 \le C_I(4K_{\hat p})^\theta N^{\theta/\hat p}\wnorm{v}{p}.
\]
Taking $p$th powers proves the result. The interpolated operator agrees with the Ritz formula on $X_p$, since both are restrictions of the same operator on the endpoint sum. In particular, $v$ need not belong to the intersection of the endpoint spaces.
\end{proof}

\begin{corollary}\label{cor:poisson-gqo}
Under Assumption~\ref{ass:discrete}, every nested sequence in $\conf(\mesh_0)$ satisfies
\begin{equation}\label{eq:poisson-gqo}
 \sum_{j=\ell}^{\ell+N-1}\wnorm{u_{j+1}-u_j}{p}^p
 \le C_{\rm gqo}N^\alpha\wnorm{u-u_\ell}{p}^p,
 \qquad u_j=R_{\mesh_j}u,
\end{equation}
for all $u\in X_p$ and $\ell,N$,
with the constants of Theorem~\ref{thm:interpolation-accumulation}.
\end{corollary}
\begin{proof}
Apply Theorem~\ref{thm:interpolation-accumulation} to $v=u-u_\ell$
and the nested sequence
$(\mesh_0,\mesh_{\ell+1},\mesh_{\ell+2},\ldots)$.
By compatibility~\eqref{eq:compatibility},
$R_{\mesh_0}v=0$ and
$R_{\mesh_{\ell+j}}v=u_{\ell+j}-u_\ell$ for $j\ge1$,
which proves the estimate.
\end{proof}

The proof also quantifies the dependence of accumulation on the discrete stability bound.

\begin{remark}[Dependence on the stability bound]
\label{rem:quantitative-accumulation}
For fixed $p\ne2$, the constants in
Theorem~\ref{thm:interpolation-accumulation} may be chosen so that
\[
 1-\alpha\ge\frac{c}{1+C_{\rm disc}},\qquad
 C_{\rm gqo}\le C(1+C_{\rm disc})^p.
\]
Here $c,C>0$ depend on $p$, $\Omega$, $r$, and the uniform mesh
and local projection bounds.
Indeed, let \(\hat p=p+{\varepsilon_p}\operatorname{sgn}(p-2)/2\)
where $\varepsilon_p$ is chosen as in the proof of
Proposition~\ref{prop:local-openness}.
The interpolation relation gives
\[
 1-\alpha=\frac{\varepsilon_p}{2|p-2|+\varepsilon_p}.
\]
The bound for $\varepsilon_p$ proves the first assertion.
As \(C_{\rm disc}\) varies, these choices of \(\hat p\) remain in a fixed compact subinterval of \((1,\infty)\).
Represent $X_2$, $X_{\hat p}$, and $X_p$ as intermediate spaces
of one fixed larger Sobolev couple.
The identifications are uniform on this compact exponent set,
so complex reiteration~\cite[Theorem~4.6.1]{BerghLofstrom1976}
bounds $C_I$ uniformly, including as $\hat p$ approaches $p$.
The estimate $K_{\hat p}\le M_p\le C(1+C_{\rm disc})$ then gives
the asserted bound for $C_{\rm gqo}$ in
\eqref{eq:interpolation-constants}.
\end{remark}

\section{Residual estimation in the target norm}
\label{sec:estimator}

We now introduce computable residual estimators and prove
global reliability and local efficiency up to data oscillation.
The resulting equivalence between the estimator and the total
error connects the quantities used for marking to the error
measured in the convergence and optimality theorems.

\subsection{Residual estimator and reliability}
For an interior face $e\in \faces{\mesh}$, write $e=\tau^+\cap\tau^-$ with $\tau^\pm\in\mesh$, 
and set $\omega_e=\tau^+\cup\tau^-$ and $h_e=|e|^{1/(d-1)}$.
Use the outward normals $n_\pm$ of $\tau^\pm$ at $e$ to define
\[
 \jump{e}(v)=\nabla v|_{\tau^+}\cdot n_+
 +\nabla v|_{\tau^-}\cdot n_-.
\]
For $v\in V_{\mesh}$ put $R_\tau(v)=(f+\Delta v)|_\tau$ and define the local error estimators by
\begin{align}
 \eta_{\mesh}(\tau,v)
 &=h_\tau^p\norm{R_\tau(v)}_{L^p(\tau)}^p
 +h_\tau\sum_{\{e\in \faces{\mesh}:e\subset\partial\tau\}}
 \norm{\jump{e}(v)}_{L^p(e)}^p,\label{eq:eta-local}\\
 \eta_{\mesh}(Q,v)&=\sum_{\tau\in Q}\eta_{\mesh}(\tau,v)
 \qquad(Q\subset \mesh).\label{eq:eta-subset}
\end{align}
Note that these are in $p$th power units.
We will write $\eta_{\mesh}(\tau)=\eta_{\mesh}(\tau,u_{\mesh})$ for the local error estimator at the Galerkin solution. 
Write also $\eta_{\mesh}(Q)=\eta_{\mesh}(Q,u_{\mesh})$ and $\eta_{\mesh}=\eta_{\mesh}(\mesh)$. Integrals over faces are with respect to surface measure. 
%There is no boundary-face residual because the homogeneous Dirichlet condition is imposed strongly.

Let $\Pi_\tau:L^p(\tau)\to\PP_{r-1}(\tau)$ be the moment projection defined by
\[
 \int_\tau(v-\Pi_\tau v)q=0\qquad\text{for all}\quad q\in\PP_{r-1}(\tau).
\]
%Although defined by the $L^2$ mass matrix, it acts on every $L^p(\tau)$, since its moments against polynomials are defined there. 
Reference element scaling gives a uniform $L^p$ bound. 
%In particular its range contains $\Delta v|_\tau$ (which is zero when $r=1$). 
Set
\begin{equation}\label{eq:osc}
 \osc_{\mesh}(\tau)=h_\tau^p\norm{f-\Pi_\tau f}_{L^p(\tau)}^p,
 \qquad
 \osc_{\mesh}(Q)=\sum_{\tau\in Q}\osc_{\mesh}(\tau).
\end{equation}
We write $ \osc_{\mesh}=\osc_{\mesh}(\mesh)$ for the whole mesh.

Elementwise integration by parts, first for smooth test functions and then by density, yields the residual representation
\begin{equation}\label{eq:residual-identity}
 a(u-u_{\mesh},v)
 =\sum_{\tau\in \mesh}\int_\tau R_\tau(u_{\mesh})v
 -\sum_{e\in \faces{\mesh}}\int_e \jump{e}(u_{\mesh})v
 \qquad(v\in X_{p'}).
\end{equation}
The face terms are well defined: $v$ has a trace, the polynomial jumps are bounded, and each face has finite measure. Galerkin orthogonality allows us to replace $v$ by $v-I_{\mesh}v$ on the right, bringing in the local approximation estimates for the projector.

\begin{theorem}[Reliability]\label{thm:reliability}
Under Assumption~\ref{ass:continuous}, we have
\begin{equation}\label{eq:reliability}
 \targeterror{p}(\mesh)^p\le C_{\rm rel}\eta_{\mesh}
 \qquad\textrm{for all}\quad \mesh\in\conf(\mesh_0).
\end{equation}
\end{theorem}
\begin{proof}
Apply \eqref{eq:continuous-stab} to $u-u_{\mesh}$ and use $a(u-u_{\mesh},w)=a(u-u_{\mesh},w-I_{\mesh}w)$ in \eqref{eq:residual-identity}.
For a volume term, Lemma~\ref{lem:local-projector} with $s=p'$ gives
\[
 \Big|\int_\tau R_\tau(w-I_{\mesh}w)\Big|
 \le C h_\tau\norm{R_\tau}_{L^p(\tau)}
 \norm{\nabla w}_{L^{p'}(\omega_{\mesh}^{m}(\tau))}.
\]
For a face term the corresponding factor is $h_e^{1/p}\norm{\jump{e}(u_\mesh)}_{L^p(e)}$, because $1-1/p'=1/p$. H\"older's inequality for the finite sums, comparability of adjacent element scales, and bounded patch overlap give
\[
 |a(u-u_{\mesh},w)|\le C\eta_{\mesh}^{1/p}\wnorm{w}{p'}.
\]
Taking the supremum and then the $p$th power proves the assertion.
\end{proof}

\subsection{Local efficiency}
Let $\omega_\tau$ be the union of $\tau$ and its face neighbors. We test the volume and face residuals using the polynomial bubbles constructed in Appendix~\ref{app:bubbles}.

\begin{theorem}[Efficiency and estimator equivalence]\label{thm:efficiency}
For every $\tau\in \mesh$, we have
\begin{equation}\label{eq:local-efficiency}
 \eta_{\mesh}(\tau)\le C_{\rm eff}
 \big(\norm{\nabla(u-u_{\mesh})}_{L^p(\omega_\tau)}^p
 +\osc_{\mesh}(\{\tau'\subset\omega_\tau\})\big).
\end{equation}
The oscillation also satisfies
\begin{equation}\label{eq:osc-domination}
 \osc_{\mesh}\le C_{\rm osc}\eta_{\mesh}.
\end{equation}
Thus under Assumption~\ref{ass:continuous}, 
the estimator $\eta_{\mesh}^{1/p}$ is equivalent to the total error
\begin{equation}\label{eq:total-error}
 \mathfrak E_p(\mesh)=
 \left(\targeterror{p}(\mesh)^p+\osc_{\mesh}\right)^{1/p}.
\end{equation}
\end{theorem}
\begin{proof}
Write $e_\mesh=u-u_{\mesh}$, $R_\tau=R_\tau(u_{\mesh})$, and
\[
 \bar R_\tau=\Pi_{\tau}f+\Delta u_{\mesh}|_\tau,
 \qquad R_\tau-\bar R_\tau=f-\Pi_{\tau}f.
\]
We now apply the bubble lemmas of Appendix~\ref{app:bubbles} with $q=p'$.
Set $w=b_\tau\bar R_\tau$, where $b_\tau$ is the product of the barycentric coordinates on $\tau$, and extend $w$ by zero outside $\tau$. By Lemma~\ref{lem:volume-bubble}, the residual identity, and H\"older's inequality, we have
\[
\begin{aligned}
 \norm{\bar R_\tau}_{L^p(\tau)}\norm{\bar R_\tau}_{L^{p'}(\tau)}
 &\lesssim \int_\tau\bar R_\tau w
 =a(e_\mesh,w)+\int_\tau(\bar R_\tau-R_\tau)w\\
 &\lesssim \big(h_\tau^{-1}\norm{\nabla e_\mesh}_{L^p(\tau)}
       +\norm{f-\Pi_\tau f}_{L^p(\tau)}\big)
       \norm{\bar R_\tau}_{L^{p'}(\tau)}.
\end{aligned}
\]
Now cancelling the last factor and using the triangle inequality, we infer
\begin{equation}\label{eq:volume-efficiency}
 h_\tau\norm{R_\tau}_{L^p(\tau)}
 \lesssim \norm{\nabla e_\mesh}_{L^p(\tau)}+\osc_{\mesh}(\tau)^{1/p} .
\end{equation}
If $\bar R_\tau=0$, this estimate follows directly from $R_\tau=f-\Pi_\tau f$.

The jump $\jump{e}=\jump{e}(u_\mesh)$ is a polynomial of degree at most $r-1$. Set $w=L_e\jump{e}$, with the lifting from Lemma~\ref{lem:face-bubble}, and extend it by zero outside $\omega_e$. Its trace vanishes on every face other than $e$. Then, similarly to the volume estimate, we obtain
\[
\begin{aligned}
 \norm{\jump{e}}_{L^p(e)}\norm{\jump{e}}_{L^{p'}(e)}
 &\lesssim\int_e\jump{e}w
 =\sum_{\tau\subset\omega_e}\int_\tau R_\tau w-a(e_\mesh,w)\\
 &\lesssim\big(h_e^{-1/p}\norm{\nabla e_\mesh}_{L^p(\omega_e)}
       +h_e^{1/p'}\sum_{\tau\subset\omega_e}\norm{R_\tau}_{L^p(\tau)}\big)
       \norm{\jump{e}}_{L^{p'}(e)}.
\end{aligned}
\]
Cancelling the last factor and using comparability of adjacent element and face scales gives
\[
 h_e^{1/p}\norm{\jump{e}}_{L^p(e)}
 \lesssim \norm{\nabla e_\mesh}_{L^p(\omega_e)}
 + \sum_{\tau\subset\omega_e}h_\tau\norm{R_\tau}_{L^p(\tau)}.
\]
The same estimate is immediate when $\jump{e}=0$. Combining it with \eqref{eq:volume-efficiency}, taking $p$th powers, and using the bounded number of faces proves \eqref{eq:local-efficiency}.

Because $\Pi_\tau\Delta u_{\mesh}=\Delta u_{\mesh}$, one has
\(f-\Pi_\tau f=(\Id-\Pi_\tau)R_\tau\).
The uniform $L^p$ bound for $\Pi_\tau$ proves \eqref{eq:osc-domination}. Summing \eqref{eq:local-efficiency} with bounded overlap and combining it with Theorem~\ref{thm:reliability} and \eqref{eq:osc-domination} proves the asserted estimator equivalence.
\end{proof}

\section{Estimator reduction and discrete reliability}
\label{sec:axioms}

With reliability and efficiency established, we turn to the change in the estimator under refinement. We separate the reduction obtained by keeping the coarse function fixed from the perturbation caused by the new Galerkin solution.
We also bound the Galerkin increment by coarse residuals near
the refined elements, providing the localization needed in
the optimality argument.

\subsection{Stability and reduction}
Let us record the useful convexity estimate
\begin{equation}\label{eq:young-power}
 |x+y|^p\le a_t|x|^p+b_t|y|^p,
\end{equation}
for $t>0$, where $a_t=(1+t)^{p-1}$ and $b_t=(1+t^{-1})^{p-1}$. 
We will apply this inequality directly to the volume and face residuals.

\begin{lemma}[Residual perturbation]\label{lem:perturbation}
For every $\mesh\in\adm(\mesh_0)$, every $v,w\in V_\mesh$,
every $Q\subset\mesh$, and every $t>0$, we have
\begin{equation}\label{eq:perturbation}
 \eta_\mesh(Q,v)\le a_t\eta_\mesh(Q,w)
       +b_t C_{\rm pert}\wnorm{v-w}{p}^p.
\end{equation}
If $\mesh'\in\adm(\mesh_0)$ and $\mesh\preceq\mesh'$,
then for every $Q\subset\mesh\cap\mesh'$ and $t>0$,
we also have
\begin{equation}\label{eq:stability}
 \eta_\mesh(Q)\le a_t\eta_{\mesh'}(Q)
 +b_t C_{\rm stab}\wnorm{u_{\mesh'}-u_{\mesh}}{p}^p.
\end{equation}
The constants are uniform in the admissible meshes and
independent of $t$.
\end{lemma}

\begin{proof}
Put $z=v-w\in V_\mesh$.
For each $\tau\in\mesh$ and each face $e\subset\partial\tau$,
polynomial inverse and trace inequalities give
\[
 h_\tau\norm{\Delta z}_{L^p(\tau)}
 \le C\norm{\nabla z}_{L^p(\tau)},\qquad
 h_\tau^{1/p}\norm{\nabla z\cdot n}_{L^p(e)}
 \le C\norm{\nabla z}_{L^p(\tau)},
\]
where $n$ is a unit normal to $e$.
Each jump is the sum of the normal traces from the two
adjacent elements, whose diameters are comparable.
Summing these bounds and using finite overlap, we obtain
\[
 \begin{aligned}
 &\sum_{\tau\in Q}h_\tau^p
       \norm{\Delta z}_{L^p(\tau)}^p
 +\sum_{\tau\in Q}h_\tau
       \sum_{e\subset\partial\tau}
       \norm{\jump{e}(z)}_{L^p(e)}^p
 \le C_{\rm pert}\wnorm{z}{p}^p.
 \end{aligned}
\]
The residual identities
$R_\tau(v)=R_\tau(w)+\Delta z$ and
$\jump{e}(v)=\jump{e}(w)+\jump{e}(z)$,
together with \eqref{eq:young-power} and the preceding bound,
prove \eqref{eq:perturbation}.

Now suppose $\mesh\preceq\mesh'$.
For a common element $\tau\in\mesh\cap\mesh'$, conformity
ensures that each face of $\tau$ remains a whole face
of the fine mesh.
The fine element across such a face lies in the corresponding
coarse neighbor. Hence a coarse function has the same volume
residual and the same flux jumps on $\tau$ when evaluated
on either mesh. Consequently, for every $v\in V_\mesh$,
we have
\[
 \eta_\mesh(Q,v)=\eta_{\mesh'}(Q,v)
 \qquad(Q\subset\mesh\cap\mesh').
\]
Since $u_\mesh\in V_{\mesh'}$, we may apply
\eqref{eq:perturbation} on the fine mesh with
$v=u_\mesh$ and $w=u_{\mesh'}$. This gives
\[
 \eta_\mesh(Q)
 =\eta_{\mesh'}(Q,u_\mesh)
 \le a_t\eta_{\mesh'}(Q)
 +b_t C_{\rm pert}\wnorm{u_{\mesh'}-u_\mesh}{p}^p.
\]
Thus \eqref{eq:stability} holds with
$C_{\rm stab}=C_{\rm pert}$.
\end{proof}

We first estimate the reduction caused by refinement for a fixed
coarse function, then use Lemma~\ref{lem:perturbation}
to account for the change in the Galerkin solution.

\begin{lemma}[Frozen and perturbed reduction]\label{lem:reduction}
Let $\mesh,\mesh'\in\adm(\mesh_0)$ with $\mesh\preceq\mesh'$,
and put $q_h=2^{-1/d}$.
For every $v\in V_{\mesh}$, we have
\begin{align}
 \eta_{\mesh'}(\mesh'\setminus\mesh,v)
 &\le q_h\eta_{\mesh}(\mesh\setminus\mesh',v),
 \label{eq:frozen-reduction}\\
 \eta_{\mesh'}(\mesh',v)
 &\le\eta_{\mesh}(\mesh,v)
 -(1-q_h)\eta_{\mesh}(\mesh\setminus\mesh',v).
 \label{eq:frozen-global}
\end{align}
For each $\delta>0$ there is $C_\delta<\infty$, independent
of the meshes, such that
\begin{equation}\label{eq:perturbed-reduction}
 \eta_{\mesh'}(\mesh'\setminus\mesh)
 \le(1+\delta)q_h\eta_{\mesh}(\mesh\setminus\mesh')
 +C_\delta\wnorm{u_{\mesh'}-u_{\mesh}}{p}^p.
\end{equation}
The coefficient $(1+\delta)q_h$ is strictly less than one
whenever $0<\delta<2^{1/d}-1$.
\end{lemma}

\begin{proof}
Fix $v\in V_{\mesh}$ and a refined coarse element
$\tau\in\mesh\setminus\mesh'$.
Its fine descendants partition $\tau$, and each satisfies
$h_{\tau'}\le q_hh_\tau$ by the volume reduction assumption.
Since the volume residual on each descendant is the
restriction of $R_\tau(v)$, we obtain
\[
 \sum_{\{\tau'\in\mesh':\tau'\subset\tau\}}
 h_{\tau'}^p\norm{R_{\tau'}(v)}_{L^p(\tau')}^p
 \le q_h^p h_\tau^p\norm{R_\tau(v)}_{L^p(\tau)}^p.
\]
Every fine face inside $\tau$ has zero jump, because $v$
is a single polynomial on $\tau$.
For an interior coarse face $e\subset\partial\tau$,
the fine faces contained in $e$ partition $e$, and their
jumps are restrictions of $\jump{e}(v)$.
The weight contributed by a descendant of $\tau$ is at most
$q_hh_\tau$. Thus the face terms contributed from this side
sum to at most
$q_hh_\tau\norm{\jump{e}(v)}_{L^p(e)}^p$.
Combining the volume and face bounds, and using
$q_h^p\le q_h$, proves \eqref{eq:frozen-reduction}
after summation over $\tau\in\mesh\setminus\mesh'$.

On common elements, the local error estimators of $v$ are unchanged.
Splitting the fine mesh into common and new elements
therefore gives
\[
 \begin{aligned}
 \eta_{\mesh'}(\mesh',v)
 &=\eta_{\mesh'}(\mesh\cap\mesh',v)
   +\eta_{\mesh'}(\mesh'\setminus\mesh,v)
 \le\eta_{\mesh}(\mesh\cap\mesh',v)
   +q_h\eta_{\mesh}(\mesh\setminus\mesh',v).
 \end{aligned}
\]
Writing $\mesh\cap\mesh'=\mesh\setminus(\mesh\setminus\mesh')$ now proves \eqref{eq:frozen-global}.

Finally, fix $\delta>0$ and choose
$t=(1+\delta)^{1/(p-1)}-1>0$, so that $a_t=1+\delta$.
Since $u_{\mesh}\in V_{\mesh'}$, the perturbation estimate
on the fine mesh gives
\[
 \eta_{\mesh'}(\mesh'\setminus\mesh)
 \le(1+\delta)\eta_{\mesh'}(\mesh'\setminus\mesh,u_{\mesh})
 +b_tC_{\rm pert}\wnorm{u_{\mesh'}-u_{\mesh}}{p}^p.
\]
Applying \eqref{eq:frozen-reduction} with $v=u_{\mesh}$
proves \eqref{eq:perturbed-reduction}, with
$C_\delta=b_tC_{\rm pert}$.
\end{proof}

\subsection{Localized discrete reliability}
We estimate the difference between nested Galerkin solutions
using coarse residuals near the refined elements.

\begin{theorem}[Discrete reliability]\label{thm:discrete-reliability}
Let $\mesh,\mesh'\in\adm(\mesh_0)$ with $\mesh\preceq\mesh'$.
Under Assumption~\ref{ass:discrete}, the set
$Q^+=Q^+(\mesh,\mesh')$ from Lemma~\ref{lem:local-projector}
satisfies
\begin{equation}\label{eq:discrete-reliability}
 \wnorm{u_{\mesh'}-u_{\mesh}}{p}^p
 \le C_{\rm drel}\eta_{\mesh}(Q^+),
\end{equation}
where $C_{\rm drel}$ is independent of the meshes.
\end{theorem}

\begin{proof}
Put $D=u_{\mesh'}-u_{\mesh}\in V_{\mesh'}$.
The discrete stability assumption on the fine mesh gives
\begin{equation}\label{eq:discrete-stab-rel-pf}
 \wnorm{D}{p}\le C_{\rm disc}
 \sup_{0\ne\hat v\in V_{\mesh'}}
 \frac{|a(D,\hat v)|}{\wnorm{\hat v}{p'}}.
\end{equation}
For $\hat v\in V_{\mesh'}$, the function
$I_{\mesh}\hat v\in V_{\mesh}$ is an admissible test function
in both Galerkin equations. Subtracting those equations gives
$a(D,I_{\mesh}\hat v)=0$.
The fine Galerkin equation therefore implies
\begin{equation}\label{eq:discrete-residual}
 a(D,\hat v)
 =F(\hat v-I_{\mesh}\hat v)
 -a(u_{\mesh},\hat v-I_{\mesh}\hat v).
\end{equation}
Set $w=\hat v-I_{\mesh}\hat v$.
By \eqref{eq:local-reproduction}, $w$ vanishes on every
coarse element outside $Q^+$.
Its trace also vanishes on each face of such an element.
Thus, after integration by parts on the coarse mesh,
the volume terms outside $Q^+$ vanish, as do all interior
face terms having an adjacent element outside $Q^+$.
Boundary terms on $\partial\Omega$ vanish because $w$
has zero trace there.

For $\tau\in Q^+$, H\"older's inequality and
\eqref{eq:sz-volume} with $s=p'$ give
\[
 \left|\int_\tau R_\tau(u_{\mesh})w\right|
 \le C h_\tau\norm{R_\tau(u_{\mesh})}_{L^p(\tau)}
 \norm{\nabla\hat v}_{L^{p'}(\omega_{\mesh}^m(\tau))}.
\]
For each remaining interior face $e$, choose an adjacent
element $\tau\in Q^+$.
Similarly, \eqref{eq:sz-face} gives
\[
 \left|\int_e\jump{e}(u_{\mesh})w\right|
 \le C h_\tau^{1/p}\norm{\jump{e}(u_{\mesh})}_{L^p(e)}
 \norm{\nabla\hat v}_{L^{p'}(\omega_{\mesh}^m(\tau))}.
\]
The $p$th powers of the residual factors are terms of
$\eta_{\mesh}(Q^+)$.
Applying H\"older's inequality to the finite sums and using
bounded patch overlap, we obtain
\[
 |a(D,\hat v)|
 \le C\eta_{\mesh}(Q^+)^{1/p}\wnorm{\hat v}{p'} ,
\]
which establishes \eqref{eq:discrete-reliability}
in light of \eqref{eq:discrete-stab-rel-pf}.
\end{proof}

Combining discrete reliability with the global frozen estimate
bounds the possible increase of the estimator under refinement.

\begin{corollary}[Quasi-monotonicity]\label{cor:monotone}
Under Assumption~\ref{ass:discrete}, there is
$C_{\rm mon}\ge1$ such that, for every
$\mesh,\mesh'\in\adm(\mesh_0)$ with $\mesh\preceq\mesh'$, we have
\begin{equation}\label{eq:quasi-monotonicity}
 \eta_{\mesh'}\le C_{\rm mon}\eta_{\mesh}.
\end{equation}
One may take
$C_{\rm mon}=2^{p-1}(1+C_{\rm pert}C_{\rm drel})$.
\end{corollary}

\begin{proof}
Put $Q^+=Q^+(\mesh,\mesh')$.
Apply \eqref{eq:perturbation} on the fine mesh with $t=1$,
$v=u_{\mesh'}$, and $w=u_{\mesh}$.
Then \eqref{eq:frozen-global} and discrete reliability give
\[
 \begin{aligned}
 \eta_{\mesh'}
 &\le 2^{p-1}\big(
       \eta_{\mesh'}(\mesh',u_{\mesh})
       +C_{\rm pert}\wnorm{u_{\mesh'}-u_{\mesh}}{p}^p
       \big)
 \le 2^{p-1}\big(
       \eta_{\mesh}
       +C_{\rm pert}C_{\rm drel}\eta_{\mesh}(Q^+)
       \big),
 \end{aligned}
\]
and $Q^+\subset\mesh$ finishes the proof.
\end{proof}

\section{The adaptive algorithm and linear convergence}
\label{sec:linear}

We now analyze the adaptive method with D\"orfler marking
in the target norm. The marking rule gives estimator reduction
up to a term involving the change in the Galerkin solution.
Combining this reduction with the accumulation estimate of
Section~\ref{sec:gqo}, we prove $R$-linear convergence of the
estimator and total error for every fixed marking parameter.

\subsection{Algorithm}
Fix $0<\theta<1$ and $C_{\rm mark}\ge1$.
Starting from $\mesh_0$, we generate meshes $\mesh_\ell$
and Galerkin solutions $u_\ell$ by the following steps.
We use the abbreviations $\eta_\ell=\eta_{\mesh_\ell}$,
$\eta_\ell(Q)=\eta_{\mesh_\ell}(Q)$, and
$\osc_\ell=\osc_{\mesh_\ell}$.

For $\ell=0,1,\ldots$, perform:
\begin{enumerate}
\item Compute $u_\ell\in V_{\mesh_\ell}$ from
\eqref{eq:galerkin}.
\item Compute the local error estimators \eqref{eq:eta-local}.
If $\eta_\ell=0$, stop.
\item Choose $M_\ell\subset\mesh_\ell$ with cardinality
minimal up to the factor $C_{\rm mark}$, subject to
\begin{equation}\label{eq:doerfler}
 \eta_\ell(M_\ell)\ge\theta\eta_\ell.
\end{equation}
\item Set
$\mesh_{\ell+1}=\operatorname{refine}(\mesh_\ell,M_\ell)$.
\end{enumerate}

The marking inequality \eqref{eq:doerfler} suffices for
linear convergence. The cardinality requirement enters
the rate optimality proof.
If the algorithm stops at level $\ell$, reliability and
\eqref{eq:osc-domination} give $u=u_\ell$ and $\osc_\ell=0$.
We then extend the sequence by keeping the mesh and solution fixed.

The marking condition ensures that refinement reduces a fixed
fraction of the estimator, up to the change in the Galerkin solution.

\begin{lemma}[One-step estimator reduction]\label{lem:one-step}
There are constants $0<\kappa<1$ and $C_{\rm red}<\infty$,
depending on $\theta$ but independent of $\ell$, such that
\begin{equation}\label{eq:one-step}
 \eta_{\ell+1}\le\kappa\eta_\ell
 +C_{\rm red}\wnorm{u_{\ell+1}-u_\ell}{p}^p
 \qquad(\ell\ge0).
\end{equation}
\end{lemma}

\begin{proof}
Every marked element is refined, so
$M_\ell\subset\mesh_\ell\setminus\mesh_{\ell+1}$.
Hence \eqref{eq:doerfler} gives
$\eta_\ell(\mesh_\ell\setminus\mesh_{\ell+1})
 \ge\theta\eta_\ell$.
Using \eqref{eq:frozen-global}, we obtain
\[
 \begin{aligned}
 \eta_{\mesh_{\ell+1}}(\mesh_{\ell+1},u_\ell)
 &\le\eta_\ell
 -(1-q_h)\eta_\ell(\mesh_\ell\setminus\mesh_{\ell+1})
 \le[1-(1-q_h)\theta]\eta_\ell.
 \end{aligned}
\]
Applying \eqref{eq:perturbation} on the fine mesh with
$v=u_{\ell+1}$ and $w=u_\ell$ therefore gives
\[
 \eta_{\ell+1}
 \le a_t[1-(1-q_h)\theta]\eta_\ell
 +b_tC_{\rm pert}\wnorm{u_{\ell+1}-u_\ell}{p}^p.
\]
Since $0<1-(1-q_h)\theta<1$ and $a_t\to1$ as
$t\downarrow0$, we can choose $t>0$ so that
$\kappa=a_t[1-(1-q_h)\theta]<1$.
Taking $C_{\rm red}=b_tC_{\rm pert}$ proves
\eqref{eq:one-step}.
\end{proof}

\subsection{Linear convergence}
We first record an elementary criterion for geometric decay
of a nonnegative sequence, following
\cite[Section~3.2]{Feischl2022}.
It combines quasi-monotonicity with a sublinear bound on
sums of consecutive terms.

\begin{lemma}[Block contraction]\label{lem:block}
Let $(x_\ell)_{\ell\ge0}$ be a nonnegative sequence.
Suppose there are $Q\ge1$ and a function
$B:\N\to(0,\infty)$ such that, for all $\ell,n\ge0$
and $N\ge1$, we have
\begin{equation}\label{eq:block-assumptions}
 x_{\ell+n}\le Qx_\ell,\qquad
 \sum_{j=\ell+1}^{\ell+N}x_j\le B(N)x_\ell.
\end{equation}
Assume that $B(N)/N\to0$ as $N\to\infty$.
Then there are constants $C<\infty$ and $q\in(0,1)$,
depending only on $Q$ and $B$, such that
\[
 x_{\ell+n}\le Cq^n x_\ell
 \qquad(\ell,n\ge0).
\]
\end{lemma}

\begin{proof}
Fix $\ell\ge0$ and $N\ge1$.
For every $j\in\{\ell+1,\ldots,\ell+N\}$,
quasi-monotonicity gives $x_{\ell+N}\le Qx_j$.
Summing these inequalities and applying the window bound,
we obtain
\[
 Nx_{\ell+N}
 \le Q\sum_{j=\ell+1}^{\ell+N}x_j
 \le QB(N)x_\ell.
\]
Since $B(N)/N\to0$, choose a fixed $N_*\ge1$ such that
$QB(N_*)/N_*\le1/2$.
Thus $x_{\ell+N_*}\le x_\ell/2$ for every $\ell$, and
iteration gives
$x_{\ell+kN_*}\le2^{-k}x_\ell$ for all $k\ge0$.

For an arbitrary $n\ge0$, write $n=kN_*+r$ with
$0\le r<N_*$.
Quasi-monotonicity controls the remaining $r$ steps, giving
\[
 x_{\ell+n}
 \le Qx_{\ell+kN_*}
 \le Q2^{-k}x_\ell
 \le 2Q\bigl(2^{-1/N_*}\bigr)^n x_\ell.
\]
We may therefore take $C=2Q$ and $q=2^{-1/N_*}$.
\end{proof}

Estimator reduction and accumulation give the window bound
for $x_\ell=\eta_\ell$, while
Corollary~\ref{cor:monotone} supplies quasi-monotonicity.

\begin{theorem}[$R$-linear convergence]\label{thm:linear}
Suppose Assumption~\ref{ass:continuous} holds.
Let $(\mesh_\ell)$ be a nested sequence satisfying the local
mesh properties of Section~\ref{sec:setting}.
Assume that \eqref{eq:discrete-stab} holds uniformly along
this sequence and that each step refines a marked set
satisfying \eqref{eq:doerfler}.
No overlay or closure estimate is required.

Then for every fixed $0<\theta<1$, there are constants
$C_{\rm lin}\ge1$ and $0<q_{\rm lin}<1$ such that
\begin{equation}\label{eq:linear}
 \eta_{\ell+n}\le C_{\rm lin}q_{\rm lin}^{\,n}\eta_\ell
 \qquad(\ell,n\ge0).
\end{equation}
After changing $C_{\rm lin}$, the same estimate holds with
$\eta_j$ replaced by $\mathfrak E_p(\mesh_j)^p$, where
$\mathfrak E_p$ is the total error defined in
\eqref{eq:total-error}.
Consequently, $\eta_\ell^{1/p}$ and
$\mathfrak E_p(\mesh_\ell)$ satisfy geometric bounds with
factor $q_{\rm lin}^{1/p}$.
In particular, $u_\ell\to u$ strongly in $X_p$ and
$\osc_\ell\to0$.
The constants inherit the dependence of $C_{\rm gqo}$
and $\alpha$ in Theorem~\ref{thm:interpolation-accumulation},
as well as the estimator constants and $\theta$.
\end{theorem}

\begin{proof}
Reindexing Corollary~\ref{cor:poisson-gqo} and using
reliability, we obtain
\begin{equation}\label{eq:estimator-gqo}
 \sum_{j=\ell+1}^{\ell+N}
 \wnorm{u_j-u_{j-1}}{p}^p
 \le C_{\rm gqo}C_{\rm rel}N^\alpha\eta_\ell,
 \qquad 0\le\alpha<1,
\end{equation}
for every $\ell\ge0$ and $N\ge1$.
For $j=\ell+1,\ldots,\ell+N$, apply
\eqref{eq:one-step} at level $j-1$.
Summing and rearranging gives
\[
 (1-\kappa)\sum_{j=\ell+1}^{\ell+N}\eta_j
 \le\kappa\eta_\ell-\kappa\eta_{\ell+N}
 +C_{\rm red}\sum_{j=\ell+1}^{\ell+N}
       \wnorm{u_j-u_{j-1}}{p}^p.
\]
Dropping the nonpositive terminal term and using
\eqref{eq:estimator-gqo}, we conclude that
\begin{equation}\label{eq:window-estimator}
 \sum_{j=\ell+1}^{\ell+N}\eta_j
 \le
 \frac{\kappa+C_{\rm red}C_{\rm gqo}C_{\rm rel}N^\alpha}
      {1-\kappa}\,\eta_\ell.
\end{equation}
Since $\alpha<1$, the factor on the right is $o(N)$.
Corollary~\ref{cor:monotone} supplies quasi-monotonicity,
so Lemma~\ref{lem:block}, applied with
$x_\ell=\eta_\ell$ and $Q=C_{\rm mon}$, proves
\eqref{eq:linear}.

The equivalence in Theorem~\ref{thm:efficiency} transfers
this estimate to $\mathfrak E_p(\mesh_\ell)^p$.
Taking $p$th roots gives the corresponding norm bounds,
and the geometric decay implies $u_\ell\to u$ in $X_p$
and $\osc_\ell\to0$.
\end{proof}

For linear convergence, the discrete stability bound need
only hold uniformly along the chosen refinement sequence.
Thus Theorem~\ref{thm:linear} also applies to a prescribed
nested sequence satisfying its mesh and marking hypotheses.
In contrast, the optimality analysis below will compare the generated
meshes with other members of the admissible family, 
and therefore the discrete stability bound must hold throughout that family.

\section{Algebraic rate optimality}
\label{sec:optimal}

Linear convergence describes decay with the number of adaptive
steps. We now measure the total error against mesh cardinality
and prove that, for sufficiently small marking parameters,
the method attains every algebraic rate available for best
approximation together with data oscillation over the
admissible family.

\subsection{Approximation classes}
We describe the algebraic rates available over the admissible
family by prescribing a budget for the number of added elements.
For $N\in\N_0$, let
\[
 \conf_N(\mesh_0)=
 \{\mesh\in\conf(\mesh_0):\#\mesh-\#\mesh_0\le N\}.
\]
For $s>0$, define
\begin{equation}\label{eq:approximation-class}
 |(u,f)|_{\mathbb A_s}
 :=\sup_{N\ge0}(N+1)^s
 \inf_{\mesh\in\conf_N(\mesh_0)}\eta_{\mesh}^{1/p}.
\end{equation}
The approximation class $\mathbb A_s$ consists of the pairs
$(u,f)$ for which this quantity is finite. Membership means
that the best estimator value under budget $N$ is bounded
by a constant times $(N+1)^{-s}$.
The term $N=0$ gives
$\eta_0^{1/p}\le |(u,f)|_{\mathbb A_s}$.

To relate this definition to best approximation in the target
norm together with data oscillation, we set
\begin{equation}\label{eq:best-total}
 \mathfrak B_p(\mesh)=
 \left(\inf_{v\in V_{\mesh}}\wnorm{u-v}{p}^p
 +\osc_{\mesh}\right)^{1/p}.
\end{equation}

\begin{proposition}[Equivalent approximation classes]
\label{prop:classes}
Under Assumptions~\ref{ass:continuous} and~\ref{ass:discrete},
the total error $\mathfrak E_p$ defined in
\eqref{eq:total-error} satisfies
\begin{equation}\label{eq:best-equivalence}
 \eta_{\mesh}^{1/p}\simeq
 \mathfrak E_p(\mesh)\simeq\mathfrak B_p(\mesh)
\end{equation}
uniformly in $\mesh$.
Replacing $\eta_{\mesh}^{1/p}$ in
\eqref{eq:approximation-class} by either
$\mathfrak E_p(\mesh)$ or $\mathfrak B_p(\mesh)$ therefore
gives the same approximation class, with equivalent
defining quantities.
\end{proposition}
\begin{proof}
The first equivalence follows from
Theorem~\ref{thm:efficiency}.
Since $u_{\mesh}\in V_{\mesh}$, we have
$\mathfrak B_p(\mesh)\le\mathfrak E_p(\mesh)$.
The reverse comparison follows from \eqref{eq:cea},
with a constant independent of $\mesh$.
Taking infima over $\conf_N(\mesh_0)$ and then the weighted
supremum over $N$ proves the assertion about approximation
classes.
\end{proof}

\subsection{Necessary marking for a strong reduction}
The approximation classes supply meshes of controlled size
and accuracy. To use them in bounding the number of marked
elements, we first show that a sufficiently strong estimator
reduction forces the enlarged refined set to satisfy the
D\"orfler condition.
In what follows, 
set
\begin{equation}\label{eq:theta-star}
 S=C_{\rm stab}C_{\rm drel},\qquad
 \theta_*=(1+S)^{-1}.
\end{equation}
Here $C_{\rm stab}$ and $C_{\rm drel}$ are the constants in
\eqref{eq:stability} and \eqref{eq:discrete-reliability},
respectively.

\begin{lemma}[Optimal-marking criterion]\label{lem:optimal-marking}
Suppose Assumption~\ref{ass:discrete} holds and fix
$0<\theta<\theta_*$.
There is $q_\theta\in(0,1)$ such that, for every
$\mesh,\mesh'\in\adm(\mesh_0)$ with $\mesh\preceq\mesh'$,
the bound
\begin{equation}\label{eq:strong-reduction}
 \eta_{\mesh'}\le q_\theta\eta_\mesh
\end{equation}
implies
\begin{equation}\label{eq:necessary-marking}
 \eta_\mesh(Q^+)\ge\theta\eta_\mesh,
 \qquad Q^+=Q^+(\mesh,\mesh').
\end{equation}
\end{lemma}
\begin{proof}
Put $Q^+=Q^+(\mesh,\mesh')$ and
$C=\mesh\setminus Q^+$.
By \eqref{eq:refined-patch}, every element of $C$ is common
to both meshes. Applying \eqref{eq:stability} on $C$ and then
\eqref{eq:discrete-reliability}, we obtain, for every $t>0$,
\begin{equation}\label{eq:marking-scalar}
 \begin{aligned}
 \eta_\mesh
 &=\eta_\mesh(Q^+)+\eta_\mesh(C)
 \le\eta_\mesh(Q^+)+a_t\eta_{\mesh'}(C)
   +b_tC_{\rm stab}\wnorm{u_{\mesh'}-u_\mesh}{p}^p\\
 &\le(1+b_tS)\eta_\mesh(Q^+)+a_t\eta_{\mesh'}.
 \end{aligned}
\end{equation}
Since $\theta<\theta_*$ and $b_t\to1$ as $t\to\infty$,
we can choose $t>0$ such that $\theta(1+b_tS)<1$.
For this choice, set
\begin{equation}\label{eq:qtheta}
 q_\theta=\frac{1-\theta(1+b_tS)}{2a_t}.
\end{equation}
Then $0<q_\theta<1$.
Using \eqref{eq:strong-reduction} in
\eqref{eq:marking-scalar}, we conclude that
\[
 \eta_\mesh(Q^+)
 \ge\frac{1-a_tq_\theta}{1+b_tS}\eta_\mesh
 \ge\theta\eta_\mesh.\qedhere
\]
\end{proof}

\subsection{Comparison meshes and complexity}
The approximation class bounds the mesh size needed to reach
a prescribed accuracy. We use a comparison mesh whose estimator
is a fixed fraction of the current value to bound the number
of marked elements. Mesh closure and linear convergence will
then turn these bounds into rate optimality for the total error.

\begin{lemma}[Cardinality of marked sets]
\label{lem:marked-cardinality}
Let $s>0$, $M=|(u,f)|_{\mathbb A_s}<\infty$, and
$0<\theta<\theta_*$.
For every adaptive level $\ell$ with $\eta_\ell>0$, we have
\begin{equation}\label{eq:marked-cardinality}
 \#M_\ell\le C M^{1/s}\eta_\ell^{-1/(sp)},
\end{equation}
where $C$ is independent of $\ell$ and may depend on
$s$, $\theta$, and $C_{\rm mark}$.
\end{lemma}
\begin{proof}
Fix a level $\ell$ with $\eta_\ell>0$ and define
$\varepsilon>0$ by
\[
 \varepsilon^p=\frac{q_\theta}{C_{\rm mon}}\eta_\ell.
\]
Quasi-monotonicity gives
$\eta_\ell\le C_{\rm mon}\eta_0\le C_{\rm mon}M^p$,
so $0<\varepsilon<M$.
By \eqref{eq:approximation-class}, we can choose
$\mesh_\varepsilon\in\adm(\mesh_0)$ such that
\[
 \eta_{\mesh_\varepsilon}\le\varepsilon^p,
 \qquad
 \#\mesh_\varepsilon-\#\mesh_0
 \le C(M/\varepsilon)^{1/s}.
\]
Let $\mesh'=\mesh_\ell\oplus\mesh_\varepsilon$.
Since $\mesh'$ refines $\mesh_\varepsilon$,
quasi-monotonicity implies
\[
 \eta_{\mesh'}
 \le C_{\rm mon}\eta_{\mesh_\varepsilon}
 \le C_{\rm mon}\varepsilon^p
 =q_\theta\eta_\ell.
\]
Lemma~\ref{lem:optimal-marking} therefore shows that
$Q^+=Q^+(\mesh_\ell,\mesh')$ satisfies the D\"orfler condition.
The nearly minimal choice of $M_\ell$, the patch bound
\eqref{eq:refined-patch}, and the overlay bound
\eqref{eq:overlay} give
\begin{align*}
 \#M_\ell
 &\le C_{\rm mark}\#Q^+
 \le C_{\rm mark}C_{\rm patch}
       \#(\mesh_\ell\setminus\mesh')
 \le C_{\rm mark}C_{\rm patch}
       (\#\mesh'-\#\mesh_\ell)\\
 &\le C_{\rm mark}C_{\rm patch}C_{\rm ov}
       (\#\mesh_\varepsilon-\#\mesh_0)
 \le C M^{1/s}\eta_\ell^{-1/(sp)}.
\end{align*}
The second line follows because each refined coarse element
is replaced by at least two descendants.
\end{proof}

Mesh closure bounds the number of added elements by the
cumulative number of marked elements. Linear convergence
allows us to control this sum in terms of the current estimator.

\begin{theorem}[Rate-optimal target-norm AFEM]\label{thm:optimal}
Assume the mesh properties of Section~\ref{sec:setting} and
Assumptions~\ref{ass:continuous} and~\ref{ass:discrete}.
If $0<\theta<\theta_*$ and marking is minimal up to the
fixed factor $C_{\rm mark}$, then for every $s>0$, we have
\begin{equation}\label{eq:rate-optimal}
 \sup_{\ell\ge0}
 (\#\mesh_\ell-\#\mesh_0+1)^{sp}\eta_\ell
 \le C_s |(u,f)|_{\mathbb A_s}^p.
\end{equation}
The same estimate holds with $\eta_\ell$ replaced by
$\mathfrak E_p(\mesh_\ell)^p$.
By Proposition~\ref{prop:classes}, the approximation class
may equivalently be defined using best $W^{1,p}$ approximation
together with data oscillation, as in \eqref{eq:best-total}.
The constant $C_s$ is independent of the data and the level.
\end{theorem}
\begin{proof}
Put $M=|(u,f)|_{\mathbb A_s}$.
It suffices to consider $M<\infty$ and a level $\ell$
with $\eta_\ell>0$.
For $0\le j<\ell$, linear convergence
\eqref{eq:linear} gives
\[
 \eta_j^{-1/(sp)}
 \le C_{\rm lin}^{1/(sp)}
 q_{\rm lin}^{(\ell-j)/(sp)}\eta_\ell^{-1/(sp)}.
\]
Mesh closure and Lemma~\ref{lem:marked-cardinality}
therefore imply
\[
 \begin{aligned}
 \#\mesh_\ell-\#\mesh_0+1
 &\le 1+C M^{1/s}
       \sum_{j=0}^{\ell-1}\eta_j^{-1/(sp)}
 \le 1+C M^{1/s}\eta_\ell^{-1/(sp)}
       \sum_{j=0}^{\ell-1}q_{\rm lin}^{(\ell-j)/(sp)}\\
 &\le C M^{1/s}\eta_\ell^{-1/(sp)}.
 \end{aligned}
\]
The last step uses the bounded geometric sum and
$\eta_\ell\le C_{\rm mon}\eta_0\le C_{\rm mon}M^p$.
This proves \eqref{eq:rate-optimal}.
At levels with $\eta_\ell=0$, the estimate is immediate.
The estimator equivalence in Theorem~\ref{thm:efficiency}
and Proposition~\ref{prop:classes} give the assertions
about the total error and the equivalent approximation class.
\end{proof}

\begin{corollary}[Convergence and optimality near $p=2$]
\label{cor:near-two-afem}
Fix a bounded polyhedral Lipschitz domain, a polynomial degree,
and a full conforming newest-vertex bisection family satisfying
the properties of Section~\ref{sec:setting}.
There is an $\varepsilon>0$, depending on the domain, degree,
and initial mesh, such that the following holds.
If $|p-2|<\varepsilon$ and $f\in L^p(\Omega)$,
the adaptive algorithm converges $R$-linearly for every fixed
$0<\theta<1$. It is algebraically rate optimal when $\theta$
is sufficiently small and the marked sets have nearly
minimal cardinality.
\end{corollary}
\begin{proof}
Combine Corollary~\ref{prop:near-two} with
Theorems~\ref{thm:linear} and~\ref{thm:optimal}.
\end{proof}

A grading condition preserves the available approximation rates
when admissible completion increases the mesh budget by at most
a fixed factor.

\begin{remark}[Graded and ungraded approximation classes]
\label{rem:graded-classes}
Let $\adm(\mesh_0)$ be a graded subfamily of an ungraded
conforming family $\mathcal U(\mesh_0)$ in the same subdivision
hierarchy. Suppose every $\mesh\in\mathcal U(\mesh_0)$ has
an admissible refinement $\hat\mesh\succeq\mesh$ satisfying
\[
 \#\hat\mesh-\#\mesh_0
 \le C_{\rm cmp}(\#\mesh-\#\mesh_0),
\]
with $C_{\rm cmp}$ independent of $\mesh$.
Then the approximation classes defined through
\eqref{eq:best-total} over the two families coincide,
with equivalent defining quantities.
Indeed, nestedness preserves every coarse approximant, while
uniform $L^p$ stability of the elementwise polynomial projections
gives $\osc_{\hat\mesh}\le C\osc_\mesh$.
Consequently,
$\mathfrak B_p(\hat\mesh)\le C\mathfrak B_p(\mesh)$.
Together with the completion bound, this transfers the
approximation rates from $\mathcal U(\mesh_0)$ to
$\adm(\mesh_0)$. The reverse comparison follows from
inclusion of the graded family in the ungraded family.

Under the stability hypotheses of Theorem~\ref{thm:optimal},
Proposition~\ref{prop:classes} therefore identifies the estimator
class on the graded family with the class defined by
\eqref{eq:best-total} over $\mathcal U(\mesh_0)$.
The discrete stability assumption is imposed on the graded family.
For compatible initial labels and a fixed grading parameter,
the graded conforming bisection families of
\cite[Appendix~A, Corollary~5]{DemlowStevenson2011}
satisfy the completion bound above.
The equivalence constants may depend on the grading parameter.
\end{remark}

\section{Numerical experiments}
\label{sec:numerics}

We apply the adaptive algorithm of Section~\ref{sec:linear} with several
choices of $p$ and measure each resulting approximation in $W^{1,q}$,
allowing $q$ to differ from $p$.
We use a problem with a reentrant corner to examine how the choice of
$p$ affects the refinement and the errors in these different norms.

\subsection{The setup}
On the L-shaped domain $\Omega_L=(-1,1)^2\setminus([0,1]\times[-1,0])$, let $(r,\varphi)$ be polar coordinates with $0<\varphi<3\pi/2$, let $\lambda=2/3$, and set
\begin{equation}\label{eq:clean-corner}
 s=r^\lambda\sin(\lambda\varphi),\qquad
 g=(1-x^2)(1-y^2),\qquad
 u=gs,\qquad
 f=-2\nabla g\cdot\nabla s-s\Delta g .
\end{equation}
The factor $s$ vanishes on the reentrant edges and $g$ on the remaining boundary, so $u\in W^{1,q}_0(\Omega_L)$ for every $q<6$. Since $s$ is harmonic and $\nabla g(0)=0$, the load is bounded with $f=O(r^{2/3})$ at the corner.

We use conforming piecewise affine elements, newest-vertex
bisection with conforming closure, homogeneous Dirichlet data, and a
sparse direct solver.
The local error estimators are computed from \eqref{eq:eta-local}.
Unless stated otherwise, we set $\theta=0.5$ in \eqref{eq:doerfler}.
We order the elements by decreasing values of the local error estimators
and mark the smallest number needed to satisfy the D\"orfler condition.
Equal values are ordered according to a fixed numbering of the elements.

We run the algorithm for
\[
 p\in\mathcal P=\{1.25,\,1.5,\,2,\,2.5,\,3,\,3.5,\,4\}.
\]
The case $p=2$ is the usual energy method.
Discrete stability (Assumption~\ref{ass:discrete}) has not been
established independently for most tested exponents, including
$p=3,3.5,4$, as discussed in Remark~\ref{rem:range}.
The value $\theta=0.5$ has also not been verified to lie below the
sufficient threshold \eqref{eq:theta-star}.
Thus our experiments at these parameter choices are exploratory.

All runs start from the same grid of width $0.25$ and continue to
$n\approx5\times10^4$.
Here $n=\dim V_{\mesh}$ is the number of degrees of freedom, equal to
the number of interior vertices.
The runs with $p=2,3,4$ are continued to $n\approx5\times10^5$.

\subsection{Error comparison}
On every mesh produced in these runs, we compute the errors
\[
 \targeterror{q}=\norm{\nabla(u-u_\mesh)}_{L^q(\Omega)}
 \qquad(q\in\mathcal P)
\]
against the exact solution \eqref{eq:clean-corner}.
For the run with exponent $p$, let $n_\ell^{[p]}$ and
$\targeterror{q,\ell}^{[p]}$ denote the degrees of freedom and the measured
error at adaptive step $\ell$, with the initial solve indexed by $\ell=0$.
We write $n_f$ for the final number of degrees of freedom in a run.

We compare the run with $p=q$ to the energy method at a budget $N$
using the ratio
\begin{equation}\label{eq:budget-envelope}
 G_q(N)=
 \frac{\min\{\targeterror{q,\ell}^{[2]}:n_\ell^{[2]}\le N\}}
      {\min\{\targeterror{q,\ell}^{[q]}:n_\ell^{[q]}\le N\}}.
\end{equation}
Thus $G_q(N)>1$ favors the run with $p=q$ in the $W^{1,q}$ error.
Taking the smallest error up to the budget allows for
nonmonotonicity along a run.

\begin{figure}[htbp]
\centering
\includegraphics[width=\linewidth]{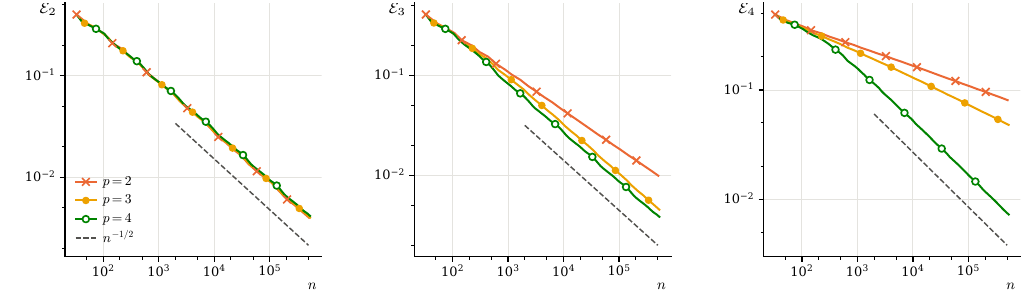}
\caption{Reentrant corner \eqref{eq:clean-corner}. Errors $\targeterror{2}$, $\targeterror{3}$, $\targeterror{4}$ from left to right, for the runs with $p=2,3,4$, with an $n^{-1/2}$ reference line.}
\label{fig:corner}
\end{figure}

Figure~\ref{fig:corner} compares the errors for the runs with $p=2,3,4$.
For $q\le2$, all tested runs show approximately $n^{-1/2}$ decay.
For the tested exponents $q\ge3$, this rate is observed only when
$p\ge q$, whereas the energy method converges more slowly as $q$
increases. This difference persists in the extended runs.
The fixed-budget gain from choosing $p=q$, shown in the left panel
of Figure~\ref{fig:gain-depth}, grows with $N$ and is much greater
for $q=4$ than for $q=3$, reaching $G_4\approx11$ at
$N=5\times10^5$.

\subsection{Refinement near the corner}
\label{sec:vertex-depth}
For \eqref{eq:clean-corner}, $u=s+(g-1)s$ with $|g-1|=O(r^2)$ and
$|\nabla g|=O(r)$, so Lemma~\ref{lem:vertex} in
Appendix~\ref{app:vertex} applies with $A=1$.
Let $h_c$ be the largest $h_\tau$ among the elements having the
reentrant corner as a vertex.
For sufficiently small $h_c$, the lemma gives a lower bound for
every discrete approximation and hence the necessary refinement
condition
\begin{equation}\label{eq:minimal-depth}
 \targeterror{4}\ge c_4h_c^{1/6},
 \qquad
 \targeterror{4}=O(n^{-1/2})\ \Longrightarrow\ h_c=O(n^{-3}).
\end{equation}

\begin{figure}[htbp]
\centering
\includegraphics[width=0.95\linewidth]{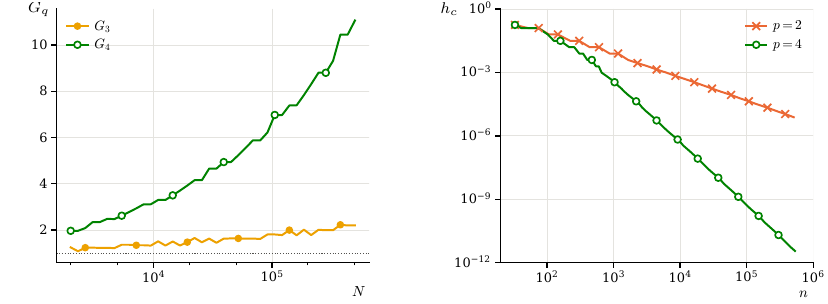}
\caption{Reentrant corner \eqref{eq:clean-corner}. Left: fixed-budget
ratios $G_3(N)$ and $G_4(N)$ of \eqref{eq:budget-envelope}. Right:
largest corner element size $h_c$ versus degrees of freedom $n$,
for the runs with $p=2,4$ at $\theta=0.5$.}
\label{fig:gain-depth}
\end{figure}

The right panel of Figure~\ref{fig:gain-depth} compares the corner
sizes for $p=2,4$ at $\theta=0.5$.
On $5\times10^4\le n\le5\times10^5$, the fitted slopes of
$\log h_c$ against $\log n$ are about $-1.1$ for $p=2$ and $-3$
for $p=4$.
In both runs, the $W^{1,4}$ error is concentrated in the corner
elements and varies approximately as $h_c^{1/6}$.
This agrees with the different error decay in Figure~\ref{fig:corner}.
The comparison is consistent with energy-norm optimality:
optimality in $W^{1,2}$ does not imply optimality in $W^{1,4}$.

\subsection{Reproducibility}
Errors are integrated against the exact solution using collapsed
tensor Gauss rules with $15\times15$ points.
Elements having the reentrant vertex are integrated by a Duffy map
collapsed at that vertex, with geometrically graded radial layers.
The singular gradient is never evaluated at the vertex.
Elements near the vertex are subdivided.

On meshes with $n\le1.6\times10^4$, raising every quadrature parameter
by one step changes $\targeterror{q}$ by at most $7\times10^{-6}$
relatively over the tested norms.
An independent quadrature implementation agrees with
$\targeterror{4}$ to $4.6\times10^{-7}$ on a mesh from the run with
$p=4$ having $h_c=5\times10^{-12}$.
Load vector and energy identity checks agree to $4.2\times10^{-10}$
and $4.9\times10^{-9}$, respectively.
The largest relative algebraic residual is approximately
$1.9\times10^{-11}$.

The implementation uses Python with NumPy and SciPy.
Every run is deterministic, and repeated and extended runs reproduce
the shorter runs bit for bit, including mesh hashes, on all common
adaptive steps.
The code used for the numerical experiments is available from the
author upon request.

\section{Conclusions and extensions}
\label{sec:conclusions}

The convergence and optimality results rely on continuous and
uniform discrete $W^{1,p}$ stability, which we establish for $p$
near $2$. Extending discrete stability to a wider range of exponents
on strongly graded adaptive meshes remains the main analytical
challenge.

For the reentrant-corner example, the lower bound shows that
half-order convergence in $W^{1,4}$ requires $h_c=O(n^{-3})$.
In the numerical experiments, $p=4$ marking approximately achieves
this refinement and the corresponding error decay, whereas
$p=2$ marking produces slower corner refinement and $W^{1,4}$
error decay over the tested range.

The analysis assumes $L^p$ data and exact Galerkin loads.
Extensions to rougher data would require suitable residual or
data-approximation estimates. Bounds on computational work would
also need to account for quadrature, inexact solves, and their costs.

\appendix
\section{Direct proof of accumulation}
\label{app:convexity}

Theorem~\ref{thm:interpolation-accumulation} was proved by interpolation. Here we give an alternative, more abstract proof based on uniform convexity of the dual Banach space. It extends the geometric argument and dyadic duality of \cite{GantumurGQO} from Hilbert spaces to this setting. Only bounded compatible projections on a single Banach space are required, and the proof below is self-contained. The adaptive argument in the main text is independent of this appendix. All spaces in this appendix are real.

Let $X$ be a Banach space, and let $\{P_k\}_{k\ge0}$ be a sequence of projections on $X$ satisfying
\[
 P_jP_k=P_{\min\{j,k\}},\qquad \sup_j\norm{P_j}\le K,\quad K\ge1.
\]
Put $P_{-1}=0$ and $\Delta_j=P_j-P_{j-1}$, so that $\norm{\Delta_j}\le2K$ and $\Delta_i\Delta_j=\delta_{ij}\Delta_j$. Suppose $X^*$ is uniformly convex and let
\[
 \delta_{X^*}(t)=\inf\big\{1-\big\|\textstyle\frac12f+\frac12g\big\|:
 \norm f,\norm g\le1,\ \norm{f-g}\ge t\big\},
 \qquad 0\le t\le2,
\]
denote its unit-ball modulus of convexity.

\begin{proposition}\label{prop:uc-accumulation}
For every $1<p<\infty$, every nonempty finite index set $I\subset\{0,1,\ldots\}$, with $n=|I|$, and every $v\in X$, we have
\begin{equation}\label{eq:uc-analysis}
 \sum_{j\in I}\norm{\Delta_jv}^p
 \le(2K)^p Bn^\beta\norm v^p,
\end{equation}
where
\[
 \mu=\min\{\delta_{X^*}(K^{-1}),\tfrac12\},\qquad
 B=1+(1-2\mu)^p<2,\qquad \beta=\log_2 B<1.
\]
\end{proposition}
\begin{proof}
Put $q=p/(p-1)$. First take $f\in\ran P_k^*$ and $g\in\ker P_k^*$. Since $P_k^*f=f$ and $P_k^*g=0$, we have
\[
 \norm f=\norm{P_k^*f}=\norm{P_k^*(f-g)}\le K\norm{f-g}.
\]
If both vectors are nonzero, the same argument applied to their normalizations gives
\[
 \Big\|\frac f{\norm f}-\frac g{\norm g}\Big\|\ge K^{-1}.
\]
The definition of the modulus and the choice of $\mu$ therefore imply
\[
 \Big\|\frac f{\norm f}+\frac g{\norm g}\Big\|\le2(1-\mu).
\]
Suppose first that $0<\norm f\le\norm g$. Writing
\[
 f+g=\norm f\Big(\frac f{\norm f}+\frac g{\norm g}\Big)
       +\Big(1-\frac{\norm f}{\norm g}\Big)g,
\]
we infer that $\norm{f+g}\le2(1-\mu)\norm f+\norm g-\norm f$. If $\norm g\le\norm f$, the analogous decomposition separates the excess of $f$ instead. Thus scalar H\"older's inequality gives the two-block estimate
\begin{equation}\label{eq:uc-two-block-short}
\begin{aligned}
 \norm{f+g}&\le\max\{\norm f,\norm g\}+(1-2\mu)\min\{\norm f,\norm g\}
 \le B^{1/p}(\norm f^q+\norm g^q)^{1/q}.
\end{aligned}
\end{equation}
This also holds if one vector vanishes, since $B\ge1$.

Order $I$ increasingly and take $z_j\in\ran\Delta_j^*$. Compatibility gives $P_k^*z_j=z_j$ for $j\le k$ and $P_k^*z_j=0$ for $j>k$. Split $I$ into two consecutive parts of almost equal size, and let $k$ be the last index in the left part. Their sums therefore lie in $\ran P_k^*$ and $\ker P_k^*$, respectively, so \eqref{eq:uc-two-block-short} applies. Repeat the split within each part. After at most $L=\lceil\log_2n\rceil$ generations, every part is a singleton. At each generation the factor is at most $B^{1/p}$ and the $q$th powers add. We conclude that
\[
 \Big\|\sum_{j\in I}z_j\Big\|\le B^{L/p}
       \Big(\sum_{j\in I}\norm{z_j}^q\Big)^{1/q}.
\]
For $n=1$, this is immediate with $L=0$.

Now put $y_j=\Delta_jv$ and $S=\sum_{j\in I}\norm{y_j}^p$. By Hahn--Banach, choose norming functionals $f_j\in X^*$ with $\norm{f_j}=\norm{y_j}^{p-1}$ and $f_j(y_j)=\norm{y_j}^p$, taking $f_j=0$ if $y_j=0$. For $z_j=\Delta_j^*f_j$, we obtain
\[
 S=\Big\langle\sum_{j\in I}z_j,v\Big\rangle,
 \qquad
 \sum_{j\in I}\norm{z_j}^q
 \le(2K)^q\sum_{j\in I}\norm{y_j}^{(p-1)q}=(2K)^q S.
\]
The synthesis bound thus yields $S\le2K B^{L/p}\norm v S^{1/q}$. If $S>0$, division by $S^{1/q}$ and taking $p$th powers give
\[
 S\le(2K)^p B^L\norm v^p\le(2K)^p Bn^\beta\norm v^p,
\]
because $L\le1+\log_2n$ and $B\ge1$. The result is immediate if $S=0$. Finally, for $I=\{\ell+1,\ldots,\ell+N\}$, compatibility gives $\Delta_j(\Id-P_\ell)u=\Delta_ju$. Applying the estimate to $v=(\Id-P_\ell)u$ gives the corresponding window bound.
\end{proof}

For $X=W^{1,p}_0(\Omega)$, the gradient is an isometry onto a closed subspace of $L^p(\Omega;\R^d)$. 
By Hahn--Banach, $X^*$ is therefore a quotient of
$L^q(\Omega;\R^d)$, where $q=p/(p-1)$.
Let $\sigma$ be normalized surface measure on $S^{d-1}$ and choose
$c_{q,d}>0$ by
$c_{q,d}^{-q}=\int_{S^{d-1}}|\omega_1|^q\,d\sigma(\omega)$.
Rotational invariance gives an isometric embedding
$F(x)\mapsto c_{q,d}F(x)\cdot\omega$
into scalar $L^q(\Omega\times S^{d-1})$.
This space is uniformly convex by Clarkson's theorem \cite{Clarkson1936},
and uniform convexity passes to closed subspaces and quotients.
For quotients, the modulus inequality applies to
norm-minimizing representatives of two unit-ball elements.
Proposition~\ref{prop:uc-accumulation} therefore gives
\eqref{eq:poisson-gqo} directly from $K_p<\infty$.

\section{Polynomial bubble functions}
\label{app:bubbles}

We provide the bubble constructions used in
Theorem~\ref{thm:efficiency}.
Throughout this appendix, $1<q<\infty$ is the exponent of
the test functions and $q'=q/(q-1)$.
The constants depend only on the polynomial degree, $q$,
the dimension, and shape regularity.

As in Section~\ref{sec:setting}, we denote by $\PP_m$ polynomials
of total degree at most $m$, on the element or face indicated
by its argument. A hat on a polynomial denotes its pullback
to the reference simplex by the affine element map.

\begin{lemma}[Volume test]\label{lem:volume-bubble}
Let $\tau$ be a shape-regular simplex, let $b_\tau$ be the
product of its barycentric coordinates, and let
$\rho\in\PP_m(\tau)$.
Then $w=b_\tau\rho\in W^{1,q}_0(\tau)$ satisfies
\begin{equation}\label{eq:volume-bubble-bounds}
 \begin{aligned}
 \int_\tau\rho w
 &\ge c\norm{\rho}_{L^{q'}(\tau)}
          \norm{\rho}_{L^q(\tau)} \qquad\text{and}\qquad
 \norm{w}_{L^q(\tau)}
 +h_\tau\norm{\nabla w}_{L^q(\tau)}
 \le C\norm{\rho}_{L^q(\tau)}.
 \end{aligned}
\end{equation}
\end{lemma}
\begin{proof}
The polynomial $w$ has zero trace.
On the fixed reference simplex $\hat\tau$, positivity
of the weighted mass form and finite dimensionality give
\begin{equation}\label{eq:weighted-polynomial}
 \int_{\hat\tau}b_{\hat\tau}\hat\rho^{\,2}
 \ge c\norm{\hat\rho}_{L^{q'}(\hat\tau)}
       \norm{\hat\rho}_{L^q(\hat\tau)}.
\end{equation}
Indeed, the quotient has a positive minimum on the compact
unit sphere of $\PP_m(\hat\tau)$.
Both sides scale by the same volume factor under an affine
map, since $1/q+1/q'=1$.
This proves the first estimate.
The second follows from boundedness of $b_\tau$ and the
polynomial inverse inequality for the fixed degree $m+d+1$.
\end{proof}

To test a face residual, we need a common polynomial trace
and a lifting that vanishes on the outer boundary of the
adjacent patch.

\begin{lemma}[Face test and lifting]\label{lem:face-bubble}
Let $e$ be an interior face, and let $b_e$ be the product
of its barycentric coordinates.
There is a linear map
$L_e:\PP_m(e)\to W^{1,q}_0(\omega_e)$ such that
$w=L_e\rho$ is a continuous piecewise polynomial with
trace $b_e\rho$ on $e$ and satisfies
\begin{equation}\label{eq:face-bubble-size}
 \begin{aligned}
 \int_e\rho w
 &\ge c\norm{\rho}_{L^{q'}(e)}
          \norm{\rho}_{L^q(e)}  \qquad\text{and}\qquad
 \norm{w}_{L^q(\omega_e)}
 +h_e\norm{\nabla w}_{L^q(\omega_e)}
 \le C h_e^{1/q}\norm{\rho}_{L^q(e)}.
 \end{aligned}
\end{equation}
\end{lemma}
\begin{proof}
Express $\rho$ as a polynomial in a fixed choice of $d-1$
barycentric coordinates of $e$.
On each adjacent simplex, extend this expression using the
corresponding barycentric coordinates and multiply by the
product of the $d$ coordinates associated with the vertices
of $e$.
The resulting polynomial equals $b_e\rho$ on $e$ and
vanishes on every other face.
The two extensions therefore agree on $e$ and have zero
trace on the outer boundary of $\omega_e$.
Their union defines $L_e\rho\in W^{1,q}_0(\omega_e)$,
and the construction is linear in $\rho$.

The first estimate follows from the trace identity and the
face version of \eqref{eq:weighted-polynomial}.
On the fixed reference face and simplex, the extension is
a bounded linear map between polynomial spaces of fixed
degree.
Affine scaling gives the factor $h_e^{1/q}$ in the
$L^q$ bound, while differentiation contributes $h_e^{-1}$.
Combining the estimates on the two adjacent simplices proves
the second bound in \eqref{eq:face-bubble-size}.
Shape regularity makes the constants uniform.
\end{proof}

\raggedbottom
\section{Corner approximation lower bound}
\label{app:vertex}

We prove the local approximation lower bound used in
Section~\ref{sec:vertex-depth} to quantify the refinement needed near
the reentrant corner.

\begin{lemma}[Vertex lower bound]\label{lem:vertex}
In two dimensions, let $V_{\mesh}$ be the conforming piecewise affine space ($r=1$). Let $c$ be a vertex of $\Omega$ with interior angle $\pi/\lambda$, $1/2<\lambda<1$, and let $s=r^\lambda\sin(\lambda\varphi)$ in polar coordinates about $c$. Suppose that $u=As+w$ on $B(c,r_0)\cap\Omega$, with $A\ne0$ and $|\nabla w(x)|\le M|x-c|^{\lambda}$ there. Let $1\le q<2/(1-\lambda)$ and $\sigma_q=\lambda-1+2/q$. There are constants $c_q>0$ and $h_0>0$, depending on $\lambda$, $q$, $A$, $M$, $r_0$ and the shape-regularity constant, such that
\[
 \norm{\nabla(u-v)}_{L^q(\tau)}\ge c_q\min\{h_\tau,h_0\}^{\sigma_q}
\]
for every $\mesh\in\conf(\mesh_0)$, every $v\in V_ \mesh$, and every $\tau\in \mesh$ having $c$ as a vertex.
\end{lemma}
\begin{proof}
Fix an element $\tau$ having $c$ as a vertex. By shape regularity,
$\operatorname{diam}\tau\le C_{\rm shape}h_\tau$. Choose
$h_0\le r_0/C_{\rm shape}$, and first suppose $h_\tau\le h_0$.
Then $\tau\subset B(c,r_0)$. Write
$\tau=c+h_\tau\hat\tau$, so that $|\hat\tau|=1$.
Homogeneity of $s$ and a change of variables give
\[
 \inf_{b\in\R^2}\norm{\nabla s-b}_{L^q(\tau)}
 =h_\tau^{\sigma_q}m(\hat\tau),
 \qquad
 m(\hat\tau)=\inf_{b\in\R^2}
 \norm{\nabla s(c+\cdot)-b}_{L^q(\hat\tau)}.
\]
The condition $q(\lambda-1)>-2$ ensures integrability at the corner.
We have $m(\hat\tau)>0$, since
$|\nabla s|=\lambda r^{\lambda-1}$ is nonconstant on every triangle.
The shape-regular triangles of unit area with a vertex at the origin
form a compact family in the closed corner sector, and $m$ depends
continuously on the triangle. Hence $m(\hat\tau)\ge m_0>0$ uniformly.

The assumed bound on $\nabla w$ gives
\[
 \norm{\nabla w}_{L^q(\tau)}
 \le CMh_\tau^{\lambda+2/q}
 =CMh_\tau^{\sigma_q+1}.
\]
Since $\nabla v$ is constant on $\tau$, the triangle inequality yields
\[
 \norm{\nabla(u-v)}_{L^q(\tau)}
 \ge |A|m_0h_\tau^{\sigma_q}
       -CMh_\tau^{\sigma_q+1}
 \ge \tfrac12|A|m_0h_\tau^{\sigma_q},
\]
after decreasing $h_0$.

If $h_\tau>h_0$, the similar subtriangle
$\tau_0=c+(h_0/h_\tau)(\tau-c)\subset\tau$
has $h_{\tau_0}=h_0$. The preceding estimate applies to $\tau_0$,
where $\nabla v$ is still constant, and therefore also gives the
required lower bound on $\tau$.
\end{proof}

\section*{Acknowledgements}
This work was partially supported by the Natural Sciences and
Engineering Research Council of Canada (NSERC) through Discovery
Grant RGPIN-2020-04389.

Claude Code was used to assist in writing the numerical code.
ChatGPT was used for language editing, literature searches, and
summaries of related work. The author checked the code, numerical
results, text, and references and takes full responsibility for the
contents of this paper.


\begin{thebibliography}{99}

\bibitem{AuscherBortzEgertSaari2019}
P.~Auscher, S.~Bortz, M.~Egert, and O.~Saari,
\emph{Nonlocal self-improving properties: a functional analytic approach},
Tunis. J. Math. \textbf{1} (2019), 151--183.

\bibitem{BadrRuss2014}
N.~Badr and E.~Russ,
\emph{A Meyers type regularity result for approximations of second order elliptic operators by $P_1$ finite elements},
Math. Nachr. \textbf{287} (2014), 498--517.
% \href{https://doi.org/10.1002/mana.201200273}{doi:10.1002/mana.201200273}.

\bibitem{BerghLofstrom1976}
J.~Bergh and J.~L\"ofstr\"om,
\emph{Interpolation Spaces: An Introduction},
Grundlehren der mathematischen Wissenschaften, vol.~223, Springer, Berlin, 1976.

\bibitem{CaoSagher1990}
W.~Cao and Y.~Sagher,
\emph{Stability of Fredholm properties on interpolation scales},
Ark. Mat. \textbf{28} (1990), 249--258.

\bibitem{CFPP2014}
C.~Carstensen, M.~Feischl, M.~Page, and D.~Praetorius,
\emph{Axioms of adaptivity},
Comput. Math. Appl. \textbf{67} (2014), 1195--1253.

\bibitem{CKNS2008}
J.~M.~Casc\'on, C.~Kreuzer, R.~H.~Nochetto, and K.~G.~Siebert,
\emph{Quasi-optimal convergence rate for an adaptive finite element method},
SIAM J. Numer. Anal. \textbf{46} (2008), 2524--2550.

\bibitem{Clarkson1936}
J.~A.~Clarkson,
\emph{Uniformly convex spaces},
Trans. Amer. Math. Soc. \textbf{40} (1936), 396--414.
% \href{https://doi.org/10.1090/S0002-9947-1936-1501880-4}{doi:10.1090/S0002-9947-1936-1501880-4}.

\bibitem{DemlowStevenson2011}
A.~Demlow and R.~Stevenson,
\emph{Convergence and quasi-optimality of an adaptive finite element method for controlling $L_2$ errors},
Numer. Math. \textbf{117} (2011), 185--218.

\bibitem{Doerfler1996}
W.~D\"orfler,
\emph{A convergent adaptive algorithm for Poisson's equation},
SIAM J. Numer. Anal. \textbf{33} (1996), 1106--1124.

\bibitem{Feischl2022}
M.~Feischl,
\emph{Inf-sup stability implies quasi-orthogonality},
Math. Comp. \textbf{91} (2022), 2059--2094.
% \href{https://doi.org/10.1090/mcom/3748}{doi:10.1090/mcom/3748}.

\bibitem{GantumurABEM2013}
T.~Gantumur,
\emph{Adaptive boundary element methods with convergence rates},
Numer. Math. \textbf{124} (2013), 471--516.

\bibitem{GantumurMultilevel2017}
T.~Gantumur,
\emph{Convergence rates of adaptive methods, Besov spaces, and multilevel approximation},
Found. Comput. Math. \textbf{17} (2017), 917--956.

\bibitem{GantumurGQO}
T.~Gantumur,
\emph{Projection geometry and relaxed quasi-orthogonality for inf-sup stable Galerkin methods},
arXiv:2609.16265, 2026.

\bibitem{JerisonKenig1995}
D.~Jerison and C.~E.~Kenig,
\emph{The inhomogeneous Dirichlet problem in Lipschitz domains},
J. Funct. Anal. \textbf{130} (1995), 161--219.

\bibitem{JinLeeLi2017}
G.~Jin, Y.~J.~Lee, and H.~Li,
\emph{Regularity estimates and optimal finite element methods in $W_p^1$ on polygonal domains},
Comput. Math. Appl. \textbf{74} (2017), 2089--2105.
% \href{https://doi.org/10.1016/j.camwa.2017.02.011}{doi:10.1016/j.camwa.2017.02.011}.

\bibitem{HLi2017}
H.~Li,
\emph{The $W_p^1$ stability of the Ritz projection on graded meshes},
Math. Comp. \textbf{86} (2017), 49--74.
% \href{https://doi.org/10.1090/mcom/3101}{doi:10.1090/mcom/3101}.

\bibitem{MNS2000}
P.~Morin, R.~H.~Nochetto, and K.~G.~Siebert,
\emph{Data oscillation and convergence of adaptive FEM},
SIAM J. Numer. Anal. \textbf{38} (2000), 466--488.

\bibitem{Sneiberg1974}
I.~Ya.~\Sneiberg,
\emph{Spectral properties of linear operators in interpolation families of Banach spaces},
Mat. Issled. \textbf{9}, no.~2(32) (1974), 214--229.

\bibitem{Stevenson2007}
R.~Stevenson,
\emph{Optimality of a standard adaptive finite element method},
Found. Comput. Math. \textbf{7} (2007), 245--269.

\bibitem{Stevenson2008}
R.~Stevenson,
\emph{The completion of locally refined simplicial partitions created by bisection},
Math. Comp. \textbf{77} (2008), 227--241.

\bibitem{YuanLi2011}
G.~Yuan and R.~Li,
\emph{Sharp a posteriori error estimate for elliptic equation with singular data},
Front. Math. China \textbf{6} (2011), 177--202.
% \href{https://doi.org/10.1007/s11464-010-0084-7}{doi:10.1007/s11464-010-0084-7}.

\end{thebibliography}
\end{document}